\documentclass[11pt,a4paper]{article}
\usepackage[T1]{fontenc}
\usepackage{lmodern}
\usepackage[a4paper,margin=26mm]{geometry}
\usepackage{amsmath,amssymb,amsthm,mathtools,mathrsfs}
\usepackage{microtype}
\usepackage{enumitem}

\usepackage{color}

\usepackage{needspace}
\usepackage[hidelinks,bookmarksnumbered]{hyperref}
\hypersetup{
 pdftitle={The one-velocity Baer--Nunziato system as a weak limit of compressible two-fluid solutions},
 pdfsubject={Homogenization from Wen's construction, conditional entropy, energy admissibility and measure-valued interpretation},
 pdfauthor={},
 pdfkeywords={two-fluid flow, effective viscous flux, Young measure, segregation, conditional entropy, pressure relaxation}
}
\numberwithin{equation}{section}
\newtheorem{theorem}{Theorem}[section]
\newtheorem{lemma}[theorem]{Lemma}
\newtheorem{proposition}[theorem]{Proposition}
\newtheorem{corollary}[theorem]{Corollary}
\theoremstyle{definition}
\newtheorem{definition}{Definition}
\theoremstyle{remark}
\newtheorem{remark}[theorem]{Remark}
\newcommand{\R}{\mathbb R}
\newcommand{\Sstress}{\mathbb S}
\newcommand{\Id}{\mathbb I}

\newcommand{\weakstar}{\stackrel{*}{\rightharpoonup}}
\newcommand{\dx}{\,\mathrm dx}

\newcommand{\dxt}{\,\mathrm dx\,\mathrm dt}

\DeclareMathOperator{\diver}{div}
\DeclareMathOperator{\tr}{tr}

\allowdisplaybreaks[2]
\DeclareMathOperator{\curl}{curl}

\newcommand{\Rbn}{\mathscr R_{\rm BN}}

\newcommand{\pbn}{p_{\rm BN}}
\title{The one-velocity Baer--Nunziato system as a weak limit of solutions to a compressible two-fluid system}
\author{{Cosmin Burtea}\thanks{{Universit\'e Paris Cit\'e and Sorbonne Universit\'e, CNRS, IMJ-PRG, F-75013 Paris, France. Email: \texttt{cosmin.burtea@imj-prg.fr}.}}
\and {Ewelina Zatorska}\thanks{{Mathematics Institute, University of Warwick, Zeeman Building, Coventry CV4 7AL, United Kingdom. Email: \texttt{ewelina.zatorska@warwick.ac.uk}.}}}
\date{17 September 2026}
\begin{document}
\maketitle
\begin{abstract}
We study how a macroscopic mixture arises from two compressible fluids
that share a velocity but occupy disjoint regions. We fix the initial
density profiles and the common initial velocity, and let the partition
between the fluids oscillate on increasingly fine scales. The fluids
remain segregated as they evolve, but in the weak limit the phase
indicator becomes a volume fraction that may lie strictly between zero
and one. For independent adiabatic exponents $\gamma^\pm\geq9/5$, we
prove that the sequence of segregated solutions converges to a global weak
solution of the viscous one-velocity Baer--Nunziato system in three
space dimensions. The mixture
pressure is the volume-weighted average of the two phase pressures,
while their difference drives the evolution of the volume fraction
along the flow and contributes to dissipation in the macroscopic energy
inequality.
\end{abstract}

\section{Introduction}\label{sec:intro}

Does a weak limit of finite-energy solutions of a compressible fluid
system still satisfy the same equations? In the classical compactness
theory for the barotropic Navier--Stokes equations, strong convergence
of the initial densities plays a role in identifying the limiting
pressure \cite{Lions1998,Hillairet2007b}. An energy bound alone allows
oscillations in the initial data, and these can persist under the
evolution. For two fluids, such oscillations arise naturally when it is
the partition between the fluids that varies on a fine scale. We ask
which equations govern the weak limit when the individual interfaces
can no longer be resolved.

This question comes from a problem of connecting two scales of modeling the compressible multi-phase flows. At a scale where each fluid can be treated as a continuum and the
interfaces can still be traced, each constituent occupies distinct
regions. We call this the \emph{mesoscopic} scale. A sharp-interface
model describes the motion in each of these regions and couples the
two fluids across their common boundary. In the absence of mass
transfer a smooth interface moves with the normal fluid velocity. The evolution of each of regions is expressed by a transport equation for the
indicator of one phase, which a-priori does not require any smoothness of the moving boundary.

In practice, the number of interfaces is or becomes too large for a model
written at this scale to be used for numerical simulations. For this reason, the description moves to a more coarser scale, where small volume may contain portions of both fluids separated by many interfaces.
 We refer to these as \emph{macroscopic models} which belong to the
family of diffuse-interface models. In this description, a single point in space records the proportion of the volume occupied by each phase which becomes the new unknown in the system, the so-called volume fraction. 
Macroscopic models are obtained in the literature either by averaging/homogenization techniques from mesoscopic models
\cite{DrewPassman1999,IshiiHibiki2011,PerrierGutierrez2021,Hillairet2007b,BreschBurteaLagoutiere2022,BreschBurteaLagoutiere2023}
or directly from Hamilton's stationary action principle
\cite{SaurelGavrilyukRenaud2003,Gavrilyuk2011,GavrilyukSaurel2002,BurteaGavrilyukPerrin2024}.

The purpose of this paper is to justify the passage from the mesoscopic to the macroscopic description in the class of finite-energy weak solutions. We fix the initial density profile of each phase and the common initial velocity, and vary only the initial partition. From the result of Wen we deduce the existence of weak solutions of two-fluid model with additive pressures and independent adiabatic exponents ($\gamma^\pm\geq9/5$), emanating from the so-called 
segregated data. We then prove that a subsequence of these weak solutions  converges to a global weak solution of the viscous one-velocity Baer–Nunziato system in three space dimensions satisfying the macroscopic energy inequality. In the limit, the phase indicators give the local volume fractions. We identify
both the pressure of the mixture and the law governing the evolution of its volume fraction.

The novelty is that this problem is for the first time analyzed in the natural setting of bounded-energy weak solutions. Strictly speaking, we investigate the limit of global-in-time weak solutions to a compressible two-fluid system, with increasingly oscillatory initial data and we identify the equations verified in the limit.

\subsection{The mesoscopic model and oscillatory data}
\label{sec:mesoscopic}
Let $\Omega\subset\R^3$ be bounded with boundary of class $C^{2+\nu}$ for
some $\nu>0$, let $T>0$, and write $Q_T=(0,T)\times\Omega$. We consider
the additive-pressure two-fluid system
\begin{equation}\label{eq:wen}
 \left\{\begin{aligned}
 \partial_t r+\diver(ru)&=0,\\
 \partial_t q+\diver(qu)&=0,\\
 \partial_t((r+q)u)+\diver((r+q)u\otimes u)
       +\nabla\bigl(p_+(r)+p_-(q)\bigr)&=\diver\Sstress(Du),
 \end{aligned}\right.
\end{equation}
with partial densities $r,q\geq0$, common velocity $u$, and no-slip boundary
condition
\begin{equation}\label{Dirichlet}
 u|_{\partial\Omega}=0.
\end{equation}
Here $Du=\tfrac12(\nabla u+\nabla u^{\mathsf T})$ and
\begin{equation}\label{eq:viscosity}
 \Sstress(Du)=2\mu Du+\lambda\tr(Du)\Id,\qquad
 \mu>0,\quad\lambda+\frac23\mu\geq0.
\end{equation}
The pressure laws are
\begin{equation}\label{eq:exponents}
 p_\pm(z)=z^{\gamma^\pm},\qquad
 \gamma^\pm\geq\frac95,\qquad z\geq0.
\end{equation}

We prescribe fixed initial profiles satisfying
\begin{equation}\label{eq:initialprofiles}
 \begin{gathered}
 \alpha^0\in L^\infty(\Omega),\qquad 0\leq\alpha^0\leq1
       \quad\text{a.e. in }\Omega,\\
 \varrho_\pm^0\in L^{\gamma^\pm}(\Omega),\qquad
 \varrho_\pm^0>0\quad\text{a.e. in }\Omega,\qquad
 u^0\in H_0^1(\Omega;\R^3).
 \end{gathered}
\end{equation}
The functions $\varrho_\pm^0$ are the initial phase densities; they may be
unbounded and need not have a uniform positive lower bound. We choose binary
partitions such that
\begin{equation}\label{eq:partitionlimit}
 \chi_k:\Omega\to\{0,1\},\qquad
 \chi_k\weakstar\alpha^0\quad\text{in }L^\infty(\Omega),
\end{equation}
and supplement \eqref{eq:wen} with
\begin{equation}\label{initial}
 \begin{aligned}
 r_k|_{t=0}&=r_k^0=\chi_k\varrho_+^0,
 &q_k|_{t=0}&=q_k^0=(1-\chi_k)\varrho_-^0,\\
 ((r_k+q_k)u_k)|_{t=0}&=m_k^0
       =(r_k^0+q_k^0)u^0.
 \end{aligned}
\end{equation}
We denote the solutions corresponding to these data by $(r_k,q_k,u_k)$. Note however that  the assumptions
\eqref{eq:initialprofiles} ensure finite energy for the initial data, and the possible weak solutions should stay in this class.

This system has been studied in \cite{Wen2021} and \cite{VasseurWenYu2019}. In particular, Wen's existence theory \cite{Wen2021} permits either phase to be absent
without requiring one partial density to dominate the other. It therefore
provides the microscopic solutions used here, together with uniform energy and higher-integrability estimates. For these solutions set
\begin{equation}\label{eq:introvariables}
 \begin{gathered}
 \varrho_k=r_k+q_k,\qquad m_k=\varrho_ku_k,\qquad
 p_k=p_+(r_k)+p_-(q_k),\\
 s_k=\begin{cases}r_k/\varrho_k,&\varrho_k>0,\\
                  0,&\varrho_k=0.
       \end{cases}
 \end{gathered}
\end{equation}
The renormalization theory for continuity and transport equations
\cite{DiPernaLions1989,NovotnyPokorny2020,Novotny2020} can be then used to derive the equation for $s_k$:
\begin{equation}\label{eq:s}
 \partial_ts_k+\diver(s_ku_k)=s_k\diver u_k,
 \qquad s_k|_{t=0}=\chi_k.
\end{equation}
What is really interesting is that at this low level of regularity, one can still show, that the initial segregation of phases $r_k^0q_k^0=0$ is preserved. Indeed, we can show that
$r_kq_k=0$, $s_k\in\{0,1\}$,
$r_k=s_k\varrho_k$ and $q_k=(1-s_k)\varrho_k$ almost everywhere. Wherever the total density is positive, the label $s_k$ specifies which fluid is present.

\medskip
Let us now observe that the weak stability question posed above has a negative answer for \eqref{eq:wen}, even for stationary solutions. Indeed, take for instance \(p_+(z)=p_-(z)=z^2\), \(\varrho_+^0=\varrho_-^0=1\) and \(u^0=0\). For each binary partition \(\chi_k\), the fields \((r_k,q_k,u_k)=(\chi_k,1-\chi_k,0)\) solve \eqref{eq:wen}: the partial densities are independent of time and the pressure \(p_k=\chi_k^2+(1-\chi_k)^2\) is identically one. Choose \(\chi_k\rightharpoonup^*\alpha^0\) as in \eqref{eq:partitionlimit}, with \(\alpha^0\) smooth, \(0<\alpha^0<1/2\) and \(\nabla\alpha^0\not\equiv0\). The weak limit is \((\alpha^0,1-\alpha^0,0)\). At zero velocity the momentum equation requires the pressure gradient to vanish, whereas the original pressure law gives
\[
\nabla\!\left[p_+(\alpha^0)+p_-(1-\alpha^0)\right]
=\nabla\!\left[(\alpha^0)^2+(1-\alpha^0)^2\right]
=(4\alpha^0-2)\nabla\alpha^0\not\equiv0.
\]
Thus the weak limit fails to satisfy \eqref{eq:wen}, although the total density is one throughout the sequence and only the partition oscillates. The density within either occupied phase is one, so averaging the phase pressures gives \(\alpha^0p_+(1)+(1-\alpha^0)p_-(1)=1\), precisely the weak limit of \(p_k\). Our main theorem resolves this problem: it identifies the limiting pressure as the volume-weighted average of the phase pressures and gives the evolution equation for the volume fraction.

\subsection{The macroscopic system and the main result}

Theorem~\ref{thm:main} extracts a subsequence for which
$(s_k,r_k,q_k,u_k)$ converges weakly to $(\alpha,R,Q,U)$ in some sense. The limiting quantities satisfy
\begin{equation}\label{eq:macro}
 \left\{\begin{aligned}
 \partial_tR+\diver(RU)&=0,\\
 \partial_tQ+\diver(QU)&=0,\\
 \partial_t\alpha+\diver(\alpha U)&=\alpha\diver U+\omega_{BN},\\
 \partial_t((R+Q)U)+\diver((R+Q)U\otimes U)+\nabla p_{BN}
       &=\diver\Sstress(DU),
 \end{aligned}\right.
\end{equation}
with zero boundary condition for $U$. Writing $U^0=u^0$, the initial data are
\begin{equation}\label{eq:introinitialmeans}
 \begin{gathered}
 R|_{t=0}=R^0=\alpha^0\varrho_+^0,\qquad
 Q|_{t=0}=Q^0=(1-\alpha^0)\varrho_-^0,\qquad
 \alpha|_{t=0}=\alpha^0,\\
 ((R+Q)U)|_{t=0}=m^0=(R^0+Q^0)U^0.
 \end{gathered}
\end{equation}
For $0<\alpha<1$, the macroscopic phase densities are
\begin{equation}\label{eq:introphasedensities}
 \varrho_+=R/\alpha,\qquad \varrho_-=Q/(1-\alpha).
\end{equation}
The pressure
and the fraction source are identified as
\begin{equation}\label{eq:introclosedp}
 \begin{aligned}
 p_{BN}&=\alpha p_+(\varrho_+)+(1-\alpha)p_-(\varrho_-),\\
 \omega_{BN}&=\frac{\alpha(1-\alpha)}{2\mu+\lambda}
                    [p_+(\varrho_+)-p_-(\varrho_-)].
 \end{aligned}
\end{equation}
These formulas explain what changes when the interfaces are no longer traceable. At the mesoscopic scale the label $s_k$ is transported by the
common velocity $u_k$, as in \eqref{eq:s}, and it only takes one of two values $s_k\in\{0,1\}$. In the macroscopic description, a point in space can accommodate both phases. This is the case whenever $\alpha\in(0,1)$. These zones represent the ``diffuse interfaces''.
The change the volume occupied by each phase is driven by the difference of the
phase pressures. System \eqref{eq:macro} is a viscous one-velocity model of
Baer--Nunziato type; the original model was introduced in
\cite{BaerNunziato1986} to describe deflagration-to-detonation transitions.

We show that the limiting quadruple $(\alpha, R,Q,U)$ is an energy-admissible weak solution to \eqref{eq:macro}. In particular, 
its energy inequality contains the
additional nonnegative dissipation
\[
 \Rbn(\alpha,R,Q)=\frac{\alpha(1-\alpha)}{2\mu+\lambda}
                 [p_+(\varrho_+)-p_-(\varrho_-)]^2.
\]
Our main result ---Theorem~\ref{thm:main} -- proves this together with the closed pressure and fraction equations, on every finite time interval.

\subsection{General background and outline of the proof}
The propagation of density oscillations in compressible Navier--Stokes
flows was studied by Serre \cite{Serre1991} in one space dimension and by
Hillairet \cite{Hillairet2007b} in several dimensions; see also
\cite{PlotnikovSokoowski2012}. When the initial densities converge
strongly, the usual compactness argument identifies the weak limit of
the pressure. For oscillatory initial data, Hillairet instead describes
this limit as the average of the pressure law against the Young measure
of the density. For regular atomic solutions with distinct density values, the weights and densities satisfy a one-velocity multifluid system.
Hillairet constructed regular solutions of that system and asked
whether such solutions could be obtained as weak limits of finite-energy
solutions of the barotropic Navier--Stokes equations.

Bresch--Huang \cite{BreschHuang2011} and Bresch--Hillairet
\cite{BreschHillairet2015} studied this question in an intermediate
regularity class, with bounded densities and
$\diver u\in L^1(0,T;L^\infty)$. In \cite{BreschHillairet2015}, the
limiting velocity and pressure are used to construct a candidate
atomic measure; uniqueness identifies it with the measure generated by
the sequence. This also removes the noncrossing condition on the phase
densities. In one space dimension the case of periodic oscillations is
related to the Bakhvalov--Eglit equations
\cite{BakhvalovEglit1983,AmosovZlotnik1996a,AmosovZlotnik1996b,Zlotnik2027}.
Further one-dimensional derivations include density-dependent
viscosities \cite{BreschHillairet2019}, distinct pressure laws
\cite{BreschBurteaLagoutiere2022,BreschBurteaLagoutiere2023}, and
thermal effects \cite{Bresch_et_al2026}. In particular, the idea of using the binary ``color'' function $s_k\in\{0,1\}$ in order to be able to prescribe different pressure was developed in the the previous three cited papers which developed on the early work of \cite{Lagoutiere2000}.

In our setting the relevant Young measure is generated by the pair
$(s_k,\varrho_k)$. The reason for keeping the phase label is already
visible in the stationary example of Section~\ref{sec:mesoscopic}.
The total density is identically one, so its distribution is just
$\delta_1$. Yet the states $(1,1)$ and $(0,1)$ occur in proportions
$\alpha^0$ and $1-\alpha^0$. The joint measure generated by $(s_k,\varrho_k)$ records this information even when the density alone cannot distinguish the phases.
Moreover, it can be identified as a convex combination of point measures. Indeed, 
the oscillatory data \eqref{initial} select one of two phase states at each point. Since $\chi_k\weakstar\alpha^0$, their initial Young measure is
\begin{equation}\label{eq:initial_mesure}
\nu_{|t=0}=\alpha^{0}\delta_{\left(  1,\varrho_{+}^{0}\right)  }+\left(
1-\alpha^{0}\right)  \delta_{\left(  0,\varrho_{-}^{0}\right)  }.
\end{equation}
Our main result shows this structure persists: for
almost every $(t,x)\in Q_T$ the joint measure is
\begin{equation}\label{structure}
 \nu_{t,x}=\alpha(t,x)\delta_{(1,\varrho_+(t,x))}
       +(1-\alpha(t,x))\delta_{(0,\varrho_-(t,x))}.
\end{equation}
Thus the phase label can still oscillate, but within each occupied
phase the density takes a single value in the limit. The weights and
the densities in \eqref{structure}, together with the common velocity
$U$, satisfy \eqref{eq:macro}--\eqref{eq:introclosedp}.
In Section~\ref{sec:younginterpretation} we interpret our findings in terms of this Young measure and write the equation satisfied by the moments of $\nu_{t,x}$, which allows to compare with the results of Hilairet from \cite{Hillairet2007b}, more directly.

Our  construction treats independent pressure laws in three dimensions using the bounds supplied by Wen's global weak solutions. \emph{It requires neither a
uniform bound on the evolving densities nor a regular macroscopic
solution. }
The main point is that we can propagate the two-state structure
of the initial Young measure in the class of bounded-energy weak
solutions. For the segregated data \eqref{eq:partitionlimit}--\eqref{initial},
this gives a positive answer to Hillairet's realization question in the
case of a common pressure law.\footnote{Hillairet \cite{Hillairet2007b}
considers the barotropic compressible Navier--Stokes equations. These
are the equations satisfied by $(\varrho_k,u_k)$ when $p_+=p_-=p$,
since segregation gives $p_k=p(\varrho_k)$.}
Indeed, the corresponding Navier--Stokes solutions generate the
two-state density law above, and the limiting phase fractions and
densities, together with the common velocity, satisfy \eqref{eq:macro} with the closure \eqref{eq:introclosedp}.
\footnote{Our proof uses the exact segregation of the approximating
data in an essential manner. It does not establish whether every sequence of weak solutions
with uniformly bounded energy whose initial data generate the same
Young measure has a subsequence whose limit satisfies the same
macroscopic equations. Indeed, prescribing the initial Young measure alone would also allow small overlaps between the phases. Our argument does not show whether such approximations lead to the same macroscopic equations.
}

 \medskip

The proof that the initial two-state structure persists for
bounded-energy weak solutions is built on two observations. Together they
replace the uniqueness argument used in \cite{BreschHillairet2015}.

\begin{enumerate}
    \item\emph{The flow preserves segregation.} It is crucial for our argument to note that  far from vacuum one has
\[
 \varrho_k s_k(1-s_k)=0\qquad\text{a.e. in }Q_T.
\]
Thus $s_k$ is either zero or one wherever $\varrho_k>0$; on vacuum we
set $s_k=0$. Where $s_k=1$, the total density is $r_k$; where $s_k=0$,
it is $q_k$. This is the formula that ensure the ``segregation'' property of the approximation, and the atomic structure of the corresponding Young measure. Consequently, for any $b\in C^{1}\left(  [0,1]\times\lbrack0,\infty)\right)  $
with $b\left(  s,0\right)  =0$, we have
\[
 b(s_k,\varrho_k)
   =s_k b(1,r_k)+(1-s_k)b(0,q_k).
\]
This identity separates the contributions of the two phases before
passage to the limit $k\to\infty$. It tells us which phase is present in the limit, but it does not rule out density oscillations within that phase.

\item\emph{The density is identified within each phase.} The second part of the proof identifies the two contributions in the
limit when $k\rightarrow\infty$ with $\alpha b\left(  1,\varrho_{+}\right)  $
and $\left(  1-\alpha\right)  b\left(  0,\varrho_{-}\right)  $ respectively.
Essentially, this identification uses a refinement of the Lions-Feireisl
effective viscous-flux argument adapted to each
of the two phases separately. As in Wen's analysis of the
two-fluid pressure \cite{Wen2021}, the difficulty is to obtain the sign of the nonlinear pressure term.
In the high
$\gamma_{\pm}$ case, we would need the following Chebyshev-type inequalities
\[
 \begin{aligned}
 \alpha\,\overline{\varrho sp_+(\varrho)}
   &\geq \overline{s\varrho}\,\overline{s p_+(\varrho)},\qquad 
 (1-\alpha)\,\overline{\varrho(1-s) p_-(\varrho)}
   &\geq \overline{(1-s)\varrho}\,\overline{(1-s)p_-(\varrho)} .
 \end{aligned}
\]
For this aspect, the segregation property proves again to be crucial.

Combined with the phasewise
effective viscous-flux identity, these inequalities show that there are no remaining density oscillations within either phase, and one has
\[
 \begin{aligned}
 s_k b(1,r_k)&\weakstar\alpha b(1,\varrho_+),\\
 (1-s_k)b(0,q_k)&\weakstar(1-\alpha)b(0,\varrho_-)
 \end{aligned}
 \qquad\text{in }L^\infty(Q_T).
\]

\item The final step of the proof is passage  to the limit in the energy inequality. To this purpose we split the viscous dissipation (at the level of phase-segregated approximation) into its curl and divergence terms. Estimating the divergence term separately in each phase gives an additional contribution. The volume-fraction equation expresses this contribution in terms of the pressure gap, giving the pressure-relaxation term in the macroscopic energy inequality.

\end{enumerate}
\medskip

We conclude introduction by comparing our result with the recent paper \cite{ChaudhuriPokornyZatorska2026}, in which a more general measure-valued notion of solution for \eqref{eq:macro}-\eqref{eq:introclosedp} was introduced. Their measure describes oscillations in the macroscopic volume
fraction, partial densities, velocity and symmetric velocity gradient.
These are the same types of variables as in our microscopic sequence
$Z_k=(s_k,r_k,q_k,u_k,Du_k)$, but the two measures arise from different
approximations. In \cite{ChaudhuriPokornyZatorska2026} the authors
approximate the macroscopic system itself: they regularize the two mass
equations and the fraction equation by artificial diffusion and use a
Faedo--Galerkin approximation for the momentum equation. They obtain
global dissipative measure-valued solutions for finite-energy initial
laws and $\gamma^\pm>1$. The two averaged mass equations hold exactly.
However, the momentum equation includes an additional concentration measure controlled by the energy defect, see \cite[Definition~1 and Theorem~3.2]{ChaudhuriPokornyZatorska2026}. 
They also showed that the Young measure can be identified with the monoatomic distribution centered on the sufficiently regular local-in-time solution, as long as the latter exists. Beyond this time, both the measure and the energy defect are  unknown quantities. 
Therefore, apart from the derivation of the system \eqref{eq:macro}, our contribution here is also strengthening the existence result from \cite{ChaudhuriPokornyZatorska2026}, at least in the case $\gamma^\pm\geq \frac95$.
 Our weak solutions also satisfy the energy inequality, and are therefore compatible with the weak-strong uniqueness principle from \cite{ChaudhuriPokornyZatorska2026}.

\paragraph{Organization of the paper.} 

Section~\ref{sec:statement} defines the weak solutions, states the
main theorem and gives its interpretation in terms of Young measures.
Section~\ref{sec:phase} explains the pressure argument first
under stronger integrability and then proves it for
$\gamma^\pm\geq9/5$; it also identifies the source in the fraction
equation. Section~\ref{sec:macroenergy} proves the energy inequality.

\section{Definitions and main result}\label{sec:statement}
We first define the pressure potential and energy of system \eqref{eq:wen}
\begin{equation}\label{eq:energydefinitions}
 H_\pm(z)=\frac{z^{\gamma^\pm}}{\gamma^\pm-1},\qquad
 \widehat E(r,q,u)=\tfrac12(r+q)|u|^2+H_+(r)+H_-(q),
\end{equation}
and
\begin{equation}
    E(r,q,u)(t)= \int_\Omega \widehat E(r,q,u)(t) \dx.
\end{equation}
A measurable triple $(r^0,q^0,m^0)$ is called finite-energy data if
$r^0,q^0\geq0$, $m^0=0$ on
$\{r^0+q^0=0\}$, and
\begin{equation}\label{eq:genericinitialenergy}
 E^0:=\int_\Omega
 \left[\frac{|m^0|^2}{2(r^0+q^0)}
                  +H_+(r^0)+H_-(q^0)\right]\dx<\infty.
\end{equation}

\subsection{Weak solutions to the two-fluid system}\label{sec:wen}
\begin{definition}[Weak solution of the additive two-fluid system]
\label{def:wen}
For finite-energy data $(r^0,q^0,m^0)$ as in
\eqref{eq:genericinitialenergy}, a measurable triple $(r,q,u)$ is a weak
solution of \eqref{eq:wen} if $r,q\geq0$ and
\begin{equation}\label{eq:basicregularity}
 \begin{gathered}
 r\in L^\infty(0,T;L^{\gamma^+}(\Omega)),\qquad
 q\in L^\infty(0,T;L^{\gamma^-}(\Omega)),\qquad {r,q\in L^2(Q_T)},\\
 \sqrt{r+q}\,u\in L^\infty(0,T;L^2),
\quad  u\in L^2(0,T;H^1_0(\Omega;\R^3)),
 \end{gathered}
\end{equation}
and the following identities hold.
\begin{enumerate}[label=\textup{(\alph*)},leftmargin=2em,itemsep=4pt]
\item \emph{Partial masses.} For $y=r,q$, with the corresponding datum
$y^0$, and every test function 
$\psi\in C_c^\infty([0,T)\times\overline\Omega)$
,
\begin{equation}\label{eq:weakmass}
 \int_{Q_T}\bigl[y\partial_t\psi+yu\cdot\nabla\psi\bigr]\dxt
       +\int_\Omega y^0\psi(0)\dx=0.
\end{equation}
In particular, the continuity equations hold after zero extension of the
densities and their fluxes outside $\Omega$.

\item \emph{Momentum.} For $m=(r+q)u$, with the corresponding datum $m^0$, and every test function $\phi\in C_c^\infty([0,T)\times\Omega;\R^3)$,
\begin{equation}\label{eq:weakmomentum}
 \begin{aligned}
 &\int_{Q_T}\bigl[m\cdot\partial_t\phi
     +(m\otimes u):\nabla\phi
     +(p_+(r)+p_-(q))\diver\phi
     -\Sstress(Du):\nabla\phi\bigr]\dxt\\
 &\hspace{30mm}+\int_\Omega m^0\cdot\phi(0)\dx=0.
 \end{aligned}
\end{equation}

\item \emph{Renormalized partial masses.} For $y=r,q$, every
$b\in C^1([0,\infty))$ with bounded derivative,
and every $\psi\in C_c^\infty([0,T)\times\overline\Omega)$,
\begin{equation}\label{eq:weakrenormalizedmass}
 \begin{aligned}
 &\int_{Q_T}\bigl[b(y)\partial_t\psi+b(y)u\cdot\nabla\psi
       +(b(y)-yb'(y))\diver u\,\psi\bigr]\dxt\\
 &\hspace{30mm}+\int_\Omega b(y^0)\psi(0)\dx=0.
 \end{aligned}
\end{equation}
\end{enumerate}
Such a weak solution is called \emph{energy-admissible} if, for almost
 every $\tau\in(0,T)$,
\begin{equation}\label{eq:wenenergy}
  E(r,q,u)(\tau)
  +\int_0^\tau\!\int_\Omega\Sstress(Du):Du\dxt\leq E^0.
\end{equation}
\end{definition}

\begin{remark}
    Definition~\ref{def:wen} covers the additive-pressure case from the work
\cite[Definition~1.1]{Wen2021} in our notation, with the distributional balances
and initial terms displayed explicitly. We also include the energy inequality
\eqref{eq:wenenergy}, retained from Wen's construction. 
\end{remark}

Below, we summarize the results from \cite[Theorem~1.2 and Lemma~5.1]{Wen2021} providing existence and higher integrability of weak solutions in the  sense of Definition~\ref{def:wen}.

\bigskip
\begin{proposition}[Wen's global weak solutions]\label{prop:wen}
Let $\Omega\subset\R^3$
be bounded with boundary of class $C^{2+\nu}$ for some $\nu>0$, let $T>0$, and $(r^0,q^0,m^0)$ be finite-energy data.
Then, there exists a global energy-admissible
weak solution in the sense of Definition~\ref{def:wen}. 
It has the following additional integrability properties
\begin{equation}\label{eq:Wengain}
 \int_{Q_T}[r^{\gamma^++\theta^+}+q^{\gamma^-+\theta^-}]\dxt\leq C_T
\end{equation}
for $\theta^+,\theta^->0$ satisfying
\begin{equation}\label{eq:gainrange}
 \theta^\pm<\gamma^\pm/3,\qquad\theta^\pm\leq\min\{1,2\gamma^\pm/3-1\}.
\end{equation}
\end{proposition}

As it is discussed for instance in \cite[Chapter 6]{NovotnyStraskraba2004}, see also \cite[Theorem 2]{NovotnyPokorny2021}, one can obtain the following extra information for the functions $\left(  r,q\right)  $
\begin{align*}
r  & \in C\left(  [0,T];L^{l_{+}}\left(  \Omega\right)  \right)  \cap
C_{weak}\left(  [0,T];L^{\gamma_{+}}\left(  \Omega\right)  \right)  ,\\
q  & \in C\left(  [0,T];L^{l_{-}}\left(  \Omega\right)  \right)  \cap
C_{weak}\left(  [0,T];L^{\gamma_{-}}\left(  \Omega\right)  \right),
\end{align*}
for any $l_{\pm}\in\lbrack1,\gamma_{\pm})$.

\subsection{Weak solution to the labeled system}\label{sec:known}
\label{sec:labelled}
We are now interested in weak solutions to the labeled system \eqref{eq:wen} with the additional transport equation \eqref{eq:s}, i.e.
\begin{equation}\label{eq:wen_s}
 \left\{\begin{aligned}
 \partial_t r+\diver(ru)&=0,\\
 \partial_t q+\diver(qu)&=0,\\
 \partial_t s+\diver(su)&=s\diver u,\\
  \partial_t((r+q)u)+\diver((r+q)u\otimes u)
       +\nabla\bigl(p_+(r)+p_-(q)\bigr)&=\diver\Sstress(Du),
 \end{aligned}\right.
\end{equation}

\begin{definition}[Weak solution of the labelled system]\label{def:labelled}
Let $(r^0,q^0,m^0)$ be finite-energy data and let
$s^0:\Omega\to\{0,1\}$ satisfy
\begin{equation}\label{initial_in}
 \begin{aligned}
r^0=s^0(r^0+q^0)\qquad \text{and}\qquad
q^0=(1-s^0)(r^0+q^0).
 \end{aligned}
\end{equation}

A quadruple $(s,r,q,u)$ is a weak solution of \eqref{eq:wen_s}, if $(r,q,u)$ is a weak solution in
Definition~\ref{def:wen}, $s\in L^\infty(Q_T)$, and, writing
$\varrho=r+q$, we have
\begin{equation}\label{eq:labelconstraints}
 s\in\{0,1\},\qquad r=s\varrho,\qquad q=(1-s)\varrho
       \quad\text{a.e. on }Q_T.
\end{equation}
Moreover, the fraction balance equation, holds in the renormalized sense
\begin{equation}\label{eq:labelledfractionweak}
 \begin{aligned}
 &\int_{Q_T}b(s)\bigl[\partial_t\psi+u\cdot\nabla\psi
                         +\diver u\,\psi\bigr]\dxt
       +\int_\Omega b(s^0)\psi(0)\dx=0
 \end{aligned}
\end{equation}
for every $b\in C^1([0,1])$ and every $\psi\in C_c^\infty([0,T)\times\overline\Omega)$.

The solution is  called \emph{energy-admissible} when it satisfies additionally \eqref{eq:wenenergy} for almost every $\tau\in(0,T)$.
\end{definition}

The existence of weak solutions in this sense follows from the Proposition \ref{prop:wen}.

\begin{corollary}[Labeled solutions for the data \eqref{initial}]\label{cor:wen}
Let $\Omega\subset\R^3$
be bounded with boundary of class $C^{2+\nu}$ for some $\nu>0$, let $T>0$, and let $(r,q,u)$ be a weak solution
from Proposition~\ref{prop:wen} with the initial data $(r^0,q^0,m^0)$, with
$m^0=(r^0+q^0)u^0$, such that there exists $s^0:\Omega\to\{0,1\}$ satisfying
\eqref{initial_in}. Then there exist a global
energy-admissible  weak solution $(s,r,q,u)$ to system \eqref{eq:wen_s} in the sense of
Definition~\ref{def:labelled}.

Moreover, for $s_k^0=\chi_k$ satisfying \eqref{eq:partitionlimit}, the corresponding sequence of weak solutions $(s_k,r_k,q_k,u_k)$ satisfies 
\begin{align}
\label{eq:energybounds}
&\|r_k\|_{L^\infty(0,T;L^{\gamma^+})}
+\|q_k\|_{L^\infty(0,T;L^{\gamma^-})}
+\|\sqrt{\varrho_k}u_k\|_{L^\infty(0,T;L^2)}
+\|u_k\|_{L^2(0,T;H^1_0)}
\leq C_T,
\\
\label{eq:extra}
&\|r_k\|_{L^2(Q_T)}
+\|q_k\|_{L^2(Q_T)}
+\|p_k\|_{L^{1+\delta}(Q_T)}
\leq C_T,
\\
\label{eq:segregation}
&r_kq_k=0,\qquad s_k\in\{0,1\},\qquad
r_k=s_k\varrho_k,\qquad q_k=(1-s_k)\varrho_k
\quad\hbox{a.e.}
\end{align}
for \[p_k=p_+(r_k)+p_-(q_k)\] and some small $\delta>0$, depending on $\theta^\pm$ defined in Proposition \ref{prop:wen}.
\end{corollary}

Derivation of equation for $s_k$ relies on the $L^2$- time-space integrability of densities $r$ and $q$ \eqref{eq:extra} and the renormalized continuity and transport equations, see  \cite[Propositions~4--5]{NovotnyPokorny2020}, \cite[Section~8.1, (8.2)--(8.8)]{Maltese_et_al2016}, and \cite[Proposition~3.3]{Novotny2020} for similar arguments. Additionally, in the latter reference,  the following regularity for $s_k$ was obtained:
\[
s_k \in C\left(  [0,T];L^{\lambda}\left(  \Omega\right)  \right)  \cap
L^{\infty}\left( 0,T;L^{\infty}\left(  \Omega\right)  \right) ,
\] 
for all $\lambda\in[1,\infty)$.
 Moreover, the product rule from \cite[Proposition~3.3]{Novotny2020} gives the equation
\[
 \partial_t\bigl(\varrho_k s_k(1-s_k)\bigr)
       +\diver\bigl(\varrho_k s_k(1-s_k)u_k\bigr)=0,
\]
satisfied in the sense of distributions.
Its initial value is zero by \eqref{initial_in}. Thus, taking a constant spatial test function one gets
\[
 \int_\Omega\varrho_k s_k(1-s_k)(t)\dx=0
\]
for all $t\in(0,T)$.
Thus $s_k\in\{0,1\}$ on
$\{\varrho_k>0\}$, and the convention on vacuum gives \eqref{eq:segregation}.




\subsection{The main result}\label{subsec:mainresult}
We first specify the constitutive functions in the weak formulation of the Baer--Nunziato system \eqref{eq:macro}.

For compatible states $0\leq\alpha\leq1$, $R,Q\geq0$, with
$R=0$ if $\alpha=0$ and $Q=0$ if $\alpha=1$, we define
\begin{equation}\label{eq:phasedensities}
 \varrho_+=\begin{cases}R/\alpha,&\alpha>0,\\0,&\alpha=0,\end{cases}
 \qquad
 \varrho_-=\begin{cases}Q/(1-\alpha),&\alpha<1,\\0,&\alpha=1.\end{cases}
\end{equation}
For such $\varrho_\pm$ we define the pressure and source terms
\begin{align}
 \pbn(\alpha,R,Q)&=\alpha p_+(\varrho_+)+(1-\alpha)p_-(\varrho_-),
       \label{eq:BNpressure}\\
 \omega_{BN}(\alpha,R,Q)&=\frac{\alpha(1-\alpha)}{2\mu+\lambda}
                  [p_+(\varrho_+)-p_-(\varrho_-)],\label{eq:BNsource}
                  \end{align}
together with pressure potential, pressure dissipation term and the total energy 
\begin{align}  
 H_{BN}(\alpha,R,Q)&=\alpha H_+(\varrho_+)+(1-\alpha)H_-(\varrho_-),
       \label{eq:BNpotential}\\
 \Rbn(\alpha,R,Q)&=\frac{\alpha(1-\alpha)}{2\mu+\lambda}
                  [p_+(\varrho_+)-p_-(\varrho_-)]^2,\label{eq:BNrelax}\\
 E_{BN}(\alpha,R,Q,U)(t)&=\int_\Omega
       [\tfrac12(R+Q)|U|^2+H_{BN}(\alpha,R,Q)](t)\dx,
       \label{eq:macroE}\\
        E_{BN}^0&=\int_\Omega\left[\frac{|m^0|^2}{2(R^0+Q^0)}
            +\alpha^0H_+(\varrho_+^0)+(1-\alpha^0)H_-(\varrho_-^0)\right]\dx.
\end{align}

\begin{definition}[Energy-admissible weak solution of the Baer--Nunziato system]
\label{def:BN}
For the fixed initial fraction $\alpha^0(x)\in[0,1]$ let the data
$(R^0,Q^0,m^0)$ be given by \eqref{eq:introinitialmeans}, with $\varrho_+^0,\varrho_-^0, m^0$ being finite-energy data in the sense of \eqref{eq:genericinitialenergy}. A quadruple
$(\alpha,R,Q,U)$ is an energy-admissible weak solution of
\eqref{eq:macro} if the following conditions hold.
\begin{enumerate}[label=\textup{(\alph*)},leftmargin=2em,itemsep=4pt]
\item\label{bn:regularity} \emph{Regularity and compatible states.}
The unknowns $(R,Q,U)$ satisfy \eqref{eq:basicregularity} with
$(r,q,u)$ replaced by $(R,Q,U)$.
Moreover $0\leq\alpha\leq1$, $R,Q\geq0$, $R=0$ on $\{\alpha=0\}$,
$Q=0$ on $\{\alpha=1\}$, and
$\pbn(\alpha,R,Q),\omega_{BN}(\alpha,R,Q)\in L^1(Q_T)$.

\item\label{bn:balances} \emph{Partial masses and momentum.}
The identities \eqref{eq:weakmass}, \eqref{eq:weakmomentum} and
\eqref{eq:weakrenormalizedmass} hold for $(R,Q,U)$, with data
$(R^0,Q^0,m^0)$, and with $\pbn(\alpha,R,Q)$ replacing
$p_+(r)+p_-(q)$ in the momentum identity.

\item\label{bn:fraction} \emph{Renormalized fraction equation.}
For every $b\in C^1([0,1])$ and every $\psi\in C_c^\infty([0,T)\times\overline\Omega)$,
\begin{equation}\label{eq:BNfractionweak}
 \begin{aligned}
 &\int_{Q_T}\bigl[b(\alpha)\partial_t\psi+b(\alpha)U\cdot\nabla\psi
       +\bigl(b(\alpha)\diver U
               +b'(\alpha)\omega_{BN}(\alpha,R,Q)\bigr)\psi\bigr]\dxt\\
 &\hspace{30mm}+\int_\Omega b(\alpha^0)\psi(0)\dx=0.
 \end{aligned}
\end{equation}

\item\label{bn:energy} \emph{Energy inequality.}
One has $H_{BN}(\alpha,R,Q)\in L^\infty(0,T;L^1)$,
$\Rbn(\alpha,R,Q)\in L^1(Q_T)$, and, for almost every $\tau\in(0,T)$,
\begin{equation}\label{eq:macroenergy}
 E_{BN}(\alpha,R,Q,U)(\tau)+\int_0^\tau\!\int_\Omega
       [\Sstress(DU):DU+\Rbn(\alpha,R,Q)]\dxt\leq E^0_{BN}.
\end{equation}
\end{enumerate}
\end{definition}

\begin{theorem}[Global weak solutions and convergence]\label{thm:main}
Let $\Omega\subset\R^3$
be bounded with boundary of class $C^{2+\nu}$ for some $\nu>0$, let $T>0$ and $(s_k,r_k,q_k,u_k)$ be the labeled weak solutions from 
Corollary~\ref{cor:wen}, with the prescribed oscillating initial data \eqref{eq:partitionlimit}, \eqref{initial}.

Then, there exist $\alpha,R,Q,U,\overline p$, such that for a subsequence still indexed by $k$, one has
\begin{equation}\label{eq:weaklimits}
 \begin{aligned}
 s_k&\weakstar\alpha &&\text{in }L^\infty(Q_T),\\
 r_k&\rightharpoonup R,\quad q_k\rightharpoonup Q
                          &&\text{in }L^2(Q_T),\\
 r_k&\weakstar R &&\text{in }L^\infty(0,T;L^{\gamma^+}(\Omega)),\\
 q_k&\weakstar Q &&\text{in }L^\infty(0,T;L^{\gamma^-}(\Omega)),\\
 u_k&\rightharpoonup U &&\text{in }L^2(0,T;H^1_0(\Omega;\R^3)),\\
 p_k&\rightharpoonup\overline p &&\text{in }L^{1+\delta}(Q_T),
 \end{aligned}
\end{equation}
where $\delta>0$ is the exponent in \eqref{eq:extra}.
The quadruple $(\alpha,R,Q,U)$ is a global energy-admissible weak solution
of the Baer--Nunziato system in the sense of Definition~\ref{def:BN}, with
initial data $(\alpha^0,R^0,Q^0,m^0)$ as in \eqref{eq:introinitialmeans}. In particular,
\begin{equation}\label{eq:closedp}
 \overline p=\pbn(\alpha,R,Q),
\end{equation}
and $\Rbn(\alpha,R,Q)\in L^1(Q_T)$.
\end{theorem}
\begin{remark}
It follows that the source term in $\alpha$ equation is square-integrable. Indeed, we have
\[
 |\omega_{BN}(\alpha,R,Q)|^2
   =\frac{\alpha(1-\alpha)}{2\mu+\lambda}\Rbn(\alpha,R,Q)
   \leq\frac1{4(2\mu+\lambda)}\Rbn(\alpha,R,Q).
\]
\end{remark}

\subsection{Instabilities, Young measures and comparison with \cite{Hillairet2007b}}\label{sec:younginterpretation}
\paragraph{Weak stability for the barotropic system.}
To compare Theorem~\ref{thm:main} with Hillairet's result
\cite{Hillairet2007b}, consider first a single pressure law
$p(z)=z^\gamma$. A weak limit of densities gives only their local
mean, and this does not in general determine the weak limit of their
pressures. Hillairet describes the missing information by the Young
measure of the density. We recall its evolution equation and then
identify this measure for the sequences constructed here.

In three space dimensions, let $\gamma>3/2$ and let
$(\varrho_k,u_k)$ be finite-energy weak solutions of
\begin{equation}\label{NavierStokes}
\left\{
\begin{array}
[c]{l}%
\partial_{t}\varrho+\operatorname{div}\left(  \varrho u\right)  =0,\\
\partial_{t}(\varrho u)+\operatorname{div}\left(  \varrho u\otimes u\right)
+\nabla\varrho^{\gamma}=\mu\Delta u+(\mu+\lambda)\nabla\operatorname{div}u,
\end{array}
\right.
\end{equation}
generated by a family of initial
data $\left(  \varrho_{k,0},u_{k,0}\right)  _{k\geq0}$ of bounded energy%
\[
 \sup_k\int_\Omega\left(\tfrac12\varrho_{k,0}|u_{k,0}|^2+
           \frac{\varrho_{k,0}^{\gamma}}{\gamma-1}\right)\dx<\infty .
\]

With strongly convergent initial densities, the classical compactness
argument identifies the limiting pressure with $p(\varrho)$
\cite{Lions1998,Hillairet2007b}. Without this preparation oscillations
may persist. A weak limit still exists along a subsequence, but its
pressure need not be $p(\varrho)$, so the limit need not satisfy
\eqref{NavierStokes}. 

This is the issue of weak stability of low-regularity solutions. For the isentropic compressible
Navier-Stokes system it has been studied first in the one-dimensional in space
case by \cite{Serre1991} and generalized in the multi-dimensional case by
\cite{Hillairet2007b} and \cite{PlotnikovSokoowski2012}. More recent results include \cite{BreschHuang2011,BreschHillairet2015,BreschHillairet2019}. 

\paragraph{Summary of Hillairet's contribution.}
One of the natural ideas applied to various nonlinear PDEs in order to construct a (low regularity and not necessarily unique) solution is to relax the notion of solution. In \cite{Hillairet2007b}, see also \cite{Serre1991} for the $1d$ case, the
author constructs a type of measure-valued solutions that take into account
only oscillations w.r.t. the density component. Let us first describe these results more precisely.

From the energy and higher-integrability bounds, after extraction of a subsequence,
\[
 \begin{aligned}
 \varrho_k&\weakstar\varrho &&\text{in }L^\infty(0,T;L^\gamma(\Omega)),\\
 u_k&\rightharpoonup u &&\text{in }L^2(0,T;H_0^1(\Omega;\mathbb R^3)),\\
 p(\varrho_k)&\rightharpoonup\overline p &&\text{in }L^{1+\delta_0}(Q_T)
 \end{aligned}
\]
for some $\delta_0>0$. We use $\nu^\varrho$ for the measure on density
values alone, to distinguish it from the joint measure $\nu$ on
phase labels and densities in \eqref{structure}.






\begin{theorem}[Hillairet {\cite[Theorem~3]{Hillairet2007b}}]
\label{Hillairet2007}
Under the preceding assumptions there is a measurable family
$\nu^\varrho_{t,x}$ of probability measures on $[0,\infty)$ such that
\begin{equation}\label{limits_presure_intro11}
  \langle\nu^\varrho,z\rangle=\varrho,\qquad
  \langle\nu^\varrho,p(z)\rangle=\overline p .
\end{equation}
Here $z$ denotes the density variable. For smooth compactly supported
$b\colon[0,\infty)\to\mathbb R$, one has in $\mathcal D'(Q_T)$
\begin{equation}\label{Kinetic111}
 \begin{aligned}
  &\partial_t\langle\nu^\varrho,b\rangle+
    \diver(\langle\nu^\varrho,b\rangle u)
    +\langle\nu^\varrho,zb'(z)-b(z)\rangle\,\diver u\\
  &\quad=\frac{1}{2\mu+\lambda}
    \big\langle\nu^\varrho,(zb'(z)-b(z))(\overline p-p(z))\big\rangle .
 \end{aligned}
\end{equation}
The limiting density and velocity satisfy, in distributions,
\begin{equation}\label{momentumpluskinetic11}
 \left\{\begin{aligned}
 \partial_t\varrho+\diver(\varrho u)&=0,\\
 \partial_t(\varrho u)+\diver(\varrho u\otimes u)+\nabla\overline p
    &=\mu\Delta u+(\mu+\lambda)\nabla\diver u .
 \end{aligned}\right.
\end{equation}
\end{theorem}

Equation \eqref{Kinetic111} encodes the oscillations of the sequence $(\varrho_k)_{k\geq 0}$. The reader will recognize that this is a reformulation, in terms of Young measures, of the renormalized continuity equation. 

Moreover, in \cite{Hillairet2007b} particular solutions for \eqref{Kinetic111}-\eqref{momentumpluskinetic11} in the form of convex combinations of Dirac masses
\begin{equation}
\nu_{  t,x  }^\varrho=\alpha\left(  t,x\right)  \delta_{\varrho
_{+}\left(  t,x\right)  }+\left(  1-\alpha\left(  t,x\right)  \right)
\delta_{\varrho_{-}(t,x)} \label{nu_intro}%
\end{equation}
were sought and it was shown that $\left(  \nu^\varrho ,u\right)  $ verifies \eqref{Kinetic111}-\eqref{momentumpluskinetic11} if and only if the functions $\left(  \alpha,\varrho_{\pm
},u\right)  $ solve the following system of PDEs:%

\begin{equation}\label{BN_one_pressure}
 \left\{\begin{aligned}
 \partial_t\alpha+u\cdot\nabla\alpha
   &=\frac{\alpha(1-\alpha)}{2\mu+\lambda}
        \bigl[p(\varrho_+)-p(\varrho_-)\bigr],\\
 \partial_t(\alpha\varrho_+)+\diver(\alpha\varrho_+u)&=0,\\
 \partial_t((1-\alpha)\varrho_-)+\diver((1-\alpha)\varrho_-u)&=0,\\
 \partial_t(\varrho u)+\diver(\varrho u\otimes u)
     +\nabla\bigl(\alpha p(\varrho_+)+(1-\alpha)p(\varrho_-)\bigr)
   &=\mu\Delta u+(\mu+\lambda)\nabla\diver u,\\
 \varrho&=\alpha\varrho_++(1-\alpha)\varrho_- .
 \end{aligned}\right.
\end{equation}

The meaning of $\left(  \text{\ref{nu_intro}}\right)  $ is that the the
sequence of densities $\left(  \varrho_{k}\right)  _{k\geq0}$
oscillates at each point $\left(  t,x\right)  $ between two states
$\varrho_{+}\left(  t,x\right)  $ and $\varrho_{-}\left(  t,x\right)  $. The distinction between the two atoms matters: if
$\varrho_+=\varrho_-$, the measure in \eqref{nu_intro} is independent
of $\alpha$ and cannot determine the fraction.

It is therefore important to  understand under which conditions the structure
of $\nu_{t,x}$ is propagated, i.e. when starting from the initial data of the form \eqref{nu_intro} (see \ref{structure}), the measure $\nu_{t,x}$ stays of the same form for latter times. This problem has been analyzed in \cite{BreschHuang2011,BreschHillairet2015} in the framework of weak solutions
with intermediate regularity introduced by D.Hoff in the $90s$ (we recall
that, roughly speaking, such solutions are generated by bounded initial
density and velocity in $H^{1}$, see \cite{Hoff1986,Hoff1995a,Hoff1995b,Hoff2002}). These solutions allow to obtain $L_{T}^{1}L_x^{\infty}$ information for the divergence of the velocity which is a key
step to be able to propagate information regarding the structure $\left(
\text{\ref{structure}}\right)  $ of the Young measures. 

The global existence of intermediate weak-solutions in the Hoff class is
ensured only for close to equlibrium configuration. 
Our research  was directly motivated by the question that was left as an open problem \cite{Hillairet2007b}: whether his solutions to system \eqref{BN_one_pressure} could be achieved as weak limits of finite energy weak solutions.

Our main result of this paper, Theorem \ref{thm:main} partially answers this question for system \eqref{eq:wen} which is more general than the one-component Nevier-Stokes system \eqref{NavierStokes}.

\paragraph{Young-measure reformulation of the main result.}

To discuss the connection between our result and the aforementioned works 
let us reformulate it in terms of Young measures generated by sequence$(s_k,\varrho_k)$. The role of the phase label $s_k$ is already visible in the stationary example above: the total density is identically one, so its distribution is just $\delta_1$. Yet the two states $(1,1)$ and $(0,1)$ occur in proportions $\alpha^0$ and $1-\alpha^0$. The joint measure records this information even when the density alone cannot distinguish the phases.
In more general setting, there exists a measurable family of probability measures on $\nu_{t,x}$ on
$[0,1]\times[0,\infty)$ (the reader is to compare with Theorem \ref{Hillairet2007}) generated by the sequence $(s_k,\varrho_k)$, such that the volume fraction, the partial
masses and the pressure have the representations
\[
 \langle\nu,s\rangle=\alpha,\qquad
 \langle\nu,sz\rangle=R,\qquad \langle\nu,(1-s)z\rangle=Q,\qquad
 \langle\nu,p_+(sz)+p_-((1-s)z)\rangle=\overline p .
\]
Above $(s,z)$ are the dummy variables for the phase label and total density.
To see how these averages evolve, we take
$b\in C_c^1([0,1]\times[0,\infty))$ and write the joint renormalized equation. Following \cite{DiPernaLions1989,Novotny2020,NovotnyPokorny2020,NovotnyPokorny2021}, we obtain
\[
\partial_{t}b\left(  s_{k},\varrho_{k}\right)  +\operatorname{div}\left(
b\left(  s_{k},\varrho_{k}\right)  u_{k}\right)  +(\varrho_{k}\frac{\partial
b}{\partial\rho}\left(  s_{k},\varrho_{k}\right)  -b\left(  s_{k},\varrho
_{k}\right)  )\operatorname{div}u_{k}=0.
\]

The product in the last term cannot in general be replaced by the
product of its weak limits. However, using the compactness of the effective flux yields in the weak limit%
\begin{multline*}
\partial_{t}\overline{b(s,\varrho)}+\operatorname{div}\left(  \overline
{b(s,\varrho)}U\right)  -\left(  \overline{\varrho\frac{\partial b}%
{\partial\rho}(s,\varrho)}-\overline{b(s,\varrho)}\right)  \operatorname{div}%
u\\
+\left(  \overline{\varrho\frac{\partial b}{\partial\rho}(s,\varrho)p\left(
s,\varrho\right)  }-\overline{b(s,\varrho)p\left(  s,\varrho\right)  }\right)
-\left(  \overline{\varrho\frac{\partial b}{\partial\rho}(s,\varrho
)}-\overline{b(s,\varrho)}\right)  \overline{p\left(  s,\varrho\right)  }=0.
\end{multline*}
Now we can interpret the weak-*
limit $\overline{b(s,\varrho)}$ as  the local average $\langle\nu_{t,x},b\rangle$, and thus the above equation becomes
\begin{equation}\label{Kinetic1}
 \begin{aligned}
 &\partial_t\langle\nu,b\rangle+\diver(\langle\nu,b\rangle U)
      +\langle\nu,z\partial_z b(s,z)-b(s,z)\rangle\,\diver U\\
 &\quad=\frac{1}{2\mu+\lambda}
   \Big\langle\nu,\bigl(z\partial_z b(s,z)-b(s,z)\bigr)
     \bigl(\overline p-p_+(sz)-p_-((1-s)z)\bigr)\Big\rangle
 \end{aligned}
\end{equation}
in $\mathcal{D}^{\prime}(Q_{T}).$
Observe that  the initial data \eqref{initial}  oscillates between two states. Indeed, for any $b\in C^{1}\left(  \left[  0,1\right]
\times\lbrack0,\infty)\right)  $ we have that%
\begin{align*}
b\left(  s_{k},\varrho_{k}\right)  _{|t=0}  & =b\left(  \chi_{k},\chi
_{k}\varrho_{+}^{0}+\left(  1-\chi_{k}\right)  \varrho_{-}^{0}\right)  \\
& =\chi_{k}b\left(  1,\varrho_{+}^{0}\right)  +\left(  1-\chi_{k}\right)
b\left(  0,\varrho_{-}^{0}\right)  .
\end{align*}
It follows that%
\[
\left\langle \nu,b\right\rangle _{|t=0}=\lim_{k\rightarrow\infty}b\left(
s_{k},\varrho_{k}\right)  _{|t=0}=\alpha^{0}b\left(  1,\varrho_{+}^{0}\right)
+\left(  1-\alpha^{0}\right)  b\left(  0,\varrho_{-}^{0}\right)  .
\]
 This means exactly that%
\begin{equation}\label{eq:initial_mesure2}
\nu_{|t=0}=\alpha^{0}\delta_{\left(  1,\varrho_{+}^{0}\right)  }+\left(
1-\alpha^{0}\right)  \delta_{\left(  0,\varrho_{-}^{0}\right)  }.
\end{equation}
Our main result shows that at positive times the same splitting is possible, since the flow
preserves segregation:
\[
 b(s_k,\varrho_k)
    =s_k b(1,\varrho_k)+(1-s_k)b(0,\varrho_k),
\]
and that, at least for a subsequence
\[
 b(s_k,\varrho_k)\weakstar
 \alpha b(1,\varrho_+)+(1-\alpha)b(0,\varrho_-)
       \qquad\text{in }L^\infty(Q_T).
\]
The right-hand side is $\langle\nu_{t,x},b\rangle$ for the measure
$\nu_{t,x}=\alpha\delta_{(1,\varrho_+)}
                 +(1-\alpha)\delta_{(0,\varrho_-)}$
in \eqref{structure}. We can now read off the nonlinear  macroscopic quantities from this measure. In particular,  the pressure on the state $(1,\rho)$
is $p_+(\rho)$, whereas on $(0,\rho)$ it is $p_-(\rho)$. Hence
\[
 \big\langle\nu,p_+(sz)+p_-((1-s)z)\big\rangle
   =\alpha p_+(\varrho_+)+(1-\alpha)p_-(\varrho_-)=\overline p .
\]
Supposing  first that $p_+=p_-=p$, and taking test functions $b(\varrho)$ depending only on the density, the derivative term in \eqref{Kinetic1} is $z b'(z)-b(z)$ and the pressure
difference is $\overline p-p(\varrho)$. We thus obtain the kinetic equation of Hillairet \eqref{Kinetic111}.

Testing with the phase label $b(s,z)=s$, instead gives
\[
 \partial_t\alpha+\diver(\alpha U)-\alpha\diver U
   =\frac{\alpha p_+(\varrho_+)-\alpha\overline p}{2\mu+\lambda}
   =\frac{\alpha(1-\alpha)}{2\mu+\lambda}
        \bigl[p_+(\varrho_+)-p_-(\varrho_-)\bigr],
\]
which is the volumetric fraction equation in the Baer-Nunziato system \eqref{eq:macro}.

\section{Proof of Theorem \ref{thm:main}--preliminaries}

\subsection{Convergence of the initial data}\label{sec:data}
Set $U^0=u^0$ and $\beta=\min\{\gamma^+,\gamma^-\}$. For the fixed
initial partial densities $\varrho_\pm^0$ in \eqref{initial}, the assumption about data being energy finite implies that
\begin{equation}\label{eq:branchinitialenergy}
 e_\pm^0:=\tfrac12\varrho_\pm^0|U^0|^2+H_\pm(\varrho_\pm^0)
 \in L^1(\Omega).
\end{equation}
The weak-* convergence \eqref{eq:partitionlimit}, then gives
\[
 \begin{aligned}
 r_k^0&\rightharpoonup R^0
      &&\text{in }L^{\gamma^+}(\Omega),\\
 q_k^0&\rightharpoonup Q^0
      &&\text{in }L^{\gamma^-}(\Omega),\\
 m_k^0&\rightharpoonup m^0=(R^0+Q^0)U^0
      &&\text{in }L^{2\beta/(\beta+1)}(\Omega).
 \end{aligned}
\]
Indeed, \eqref{eq:branchinitialenergy} gives
$\varrho_\pm^0\in L^{\gamma^\pm}$ and
$\varrho_\pm^0U^0\in L^{2\gamma^\pm/(\gamma^\pm+1)}$.
Noticing that the initial energy depends on $k$ only linearly w.r.t. $\chi_k$ and using $\chi_k\weakstar\alpha^0$ we therefore deduce
\begin{equation}\label{eq:initialenergyconv}
 E_k^0=\int_\Omega\bigl[\chi_ke_+^0+(1-\chi_k)e_-^0\bigr]\dx
 \longrightarrow
 \int_\Omega\bigl[\alpha^0e_+^0+(1-\alpha^0)e_-^0\bigr]\dx=E^0_{BN}.
\end{equation}
In particular, the initial energies are uniformly bounded, as required
for the estimates in Corollary~\ref{cor:wen}.

\subsection{Convergence following from the uniform estimates}
\label{subsec:preliminarycompactness}
\label{sec:compactness}
Recall the notation
\[
 \varrho_k=r_k+q_k,\qquad m_k=\varrho_ku_k,\qquad
 p_k=p_+(r_k)+p_-(q_k).
\]
We may now use the uniform estimates \eqref{eq:energybounds} and
\eqref{eq:extra} from Corollary~\ref{cor:wen} to extract a subsequence
such that all convergences specified in \eqref{eq:weaklimits} hold.
We use a standard compactness mechanism to identify the limits of the
transport fluxes; see, for instance, \cite{FeireislPokorny}.

With $\beta$ as above, define
\begin{equation}\label{eq:compactexponents}
 p_0=\frac{2\beta}{\beta+1}>\frac65,\qquad
 q_0=\frac{6\beta}{4\beta+3}>1.
\end{equation}
H\"older's inequality implies that
\begin{equation}\label{eq:momentumbounds}
 m_k\text{ bounded in }L^\infty(0,T;L^{p_0}),\qquad
 m_k\otimes u_k\text{ bounded in }L^2(0,T;L^{q_0}),
\end{equation}
uniformly in $k$. 
Moreover,  $p_k$ is bounded in $L^{1+\delta}(Q_T)$, thus using 
\eqref{eq:momentumbounds} and the momentum equation, we can show that $\partial_tm_k$ is bounded
locally in $L^{\min\{2,1+\delta\}}(0,T;H^{-\ell})$ for a sufficiently
large integer $\ell$, which gives equicontinuity of $\int_{\Omega}m_k\varphi\ \dx$, for any $\varphi\in C^\infty_c(\Omega)$. The Sobolev imbedding then implies
\begin{equation}\label{eq:mstrong}
 m_k\longrightarrow m\quad\text{in }
 C([0,T];H^{-1}_{\mathrm{loc}}(\Omega)).
\end{equation}
The same argument for each continuity equation gives
$r_k\to R$, $q_k\to Q$ and $\varrho_k\to R+Q$ in $ C([0,T];H^{-1}_{\mathrm{loc}}(\Omega))$.
Pairing these convergences with the weak convergence of $u_k$ in
$L^2H^1$ identifies
\begin{equation}\label{eq:convectiveclosure}
\begin{aligned}
r_ku_k&\rightharpoonup RU
&&\text{in }L^2\!\left(0,T;L^{\frac{6\gamma^+}{\gamma^++6}}(\Omega;\mathbb R^3)\right),\\
q_ku_k&\rightharpoonup QU
&&\text{in }L^2\!\left(0,T;L^{\frac{6\gamma^-}{\gamma^-+6}}(\Omega;\mathbb R^3)\right),\\
\varrho_k u_k&\overset{*}{\rightharpoonup}(R+Q)U
&&\text{in }L^\infty(0,T;L^{p_0}(\Omega;\mathbb R^3)),\\
\varrho_k u_k\otimes u_k&\rightharpoonup(R+Q)U\otimes U
&&\text{in }L^2(0,T;L^{q_0}(\Omega;\mathbb R^{3\times3})),
\end{aligned}
\end{equation}
where \(p_0,\ q_0\) are defined in \eqref{eq:compactexponents}.
Similarly, we identify
\begin{equation}\label{eq:fractionvelocity}
 s_ku_k\rightharpoonup\alpha U\quad\text{in }L^2(Q_T).
\end{equation}

We can also use the time-derivative estimates, together with the uniform $L^\infty(L^p(\Omega)$, to extend the above convergences to:
\begin{equation}\label{conv_faible}
 \begin{aligned}
 \varrho_k u_k&\longrightarrow (R+Q)U
    &&\text{in }C_{\rm weak}([0,T];L^{p_0}(\Omega)),\\
 r_k&\longrightarrow R
    &&\text{in }C_{\rm weak}([0,T];L^{\gamma^+}(\Omega)),\\
 q_k&\longrightarrow Q
    &&\text{in }C_{\rm weak}([0,T];L^{\gamma^-}(\Omega)),\\
 s_k&\longrightarrow\alpha
    &&\text{in }C_{\rm weak}([0,T];L^p(\Omega)),\quad 1<p<\infty.
 \end{aligned}
\end{equation}

We also have
\begin{equation}\label{lim_pres}
    p_k\rightharpoonup \overline p= \overline{p_+(r)}+\overline{p_-(q)}\quad\text{in }L^{1+\delta}(Q_T), 
\end{equation}
and
\begin{equation}\label{eq:G}
 s_k\diver u_k\rightharpoonup\overline{s \diver u}:=G\quad\text{in }L^2(Q_T).
\end{equation}

Passage to the limit in the masses and volumetric fraction therefore gives
\begin{equation}\label{eq:limRQa}
 \begin{aligned}
  \int_{Q_T}\bigl[R\partial_t\psi+R U\cdot\nabla\psi\bigr]\dxt+\int_\Omega R^0\psi(0)\dx&=0,\\
    \int_{Q_T}\bigl[Q\partial_t\psi+Q U\cdot\nabla\psi\bigr]\dxt
       +\int_\Omega Q^0\psi(0)\dx&=0,\\
          \int_{Q_T}\bigl[\alpha\partial_t\psi+\alpha U\cdot\nabla\psi\bigr]\dxt
       +\int_\Omega \alpha^0\psi(0)\dx&=-\int_{Q_T}G\psi\dxt,
 \end{aligned}
\end{equation}
for any $\psi\in C_c^\infty([0,T)\times\overline\Omega)$, and in the momentum equation 
\begin{equation}\label{eq:limmomentum}
 \begin{aligned}
 &\int_{Q_T}\bigl[(R+Q)U\cdot\partial_t\phi
     +((R+Q)U\otimes U):\nabla\phi
     +\overline p\diver\phi
     -\Sstress(DU):\nabla\phi\bigr]\dxt\\
 &\hspace{30mm}+\int_\Omega m^0\cdot\phi(0)\dx=0,
 \end{aligned}
\end{equation}
for any $\phi\in C_c^\infty([0,T)\times\Omega;\R^3)$.

Finally, extending $R,Q,\alpha,U, G$ by zero outside $\Omega$ and using \cite[Lemma 6.8]{NovotnyStraskraba2004}, we can show that equations for $R,Q,$ and $\alpha$ are satisfied in the sense of distributions on $(0,T)\times \R^3$.

\subsection{Definition and estimates for partial densities \texorpdfstring{$\varrho_\pm$}{varrho±}}
We start this section from checking that  
\begin{equation}\label{cond}
    R=0\ \text{on }\{\alpha=0\},\qquad
Q=0\ \text{on }\{\alpha=1\}.
\end{equation}
Indeed, a more general result holds : consider $b:[0,\infty)\rightarrow[0,\infty)$ with $b(0)=0$ and assume there is some $l\in(1,\infty)$ such that
\begin{equation}
\label{eq:compatibility_func_of_r_1}
b(r_k)\rightharpoonup\overline{b(r)}\text{ weakly in } L^{\frac{l}{l-1}}(Q_T).
\end{equation}
Then 
\begin{equation}\label{eq:compatibility_func_of_r_2}
    \overline{b(r)}=0 \text{ a.e. on } \{\alpha=0\}.
\end{equation}
Indeed, owing to separation of phases we have \(b(r_k)=s_kb(r_k)\) and for all $l\in[1,\infty)$, \(s_k^l=s_k\), one can estimate
\[
\int_{\{\alpha=0\}} b(r_k)
=\int_{\{\alpha=0\}} s_kb(r_k)
\leq \|b(r_k)\|_{L^{\frac{l}{l-1}}(Q_T)}
       \left(\int_{\{\alpha=0\}} s_k\right)^{1/l}
\longrightarrow0,
\]
because \(s_k\rightharpoonup^\ast\alpha\) in $L^\infty(Q_T)$, and \(\int_{\{\alpha=0\}}\alpha=0\). Passing to the weak limit gives \(\int_{\{\alpha=0\}} \overline{b(r)}=0\), hence \(\overline{b(r)}=0\) a.e. on \({\{\alpha=0\}}\). The argument for \(Q\) is identical.

Conditions \eqref{cond} allow us to define \(\varrho_\pm\) as in \(\eqref{eq:phasedensities}\), namely
\begin{equation}\label{eq:phasedensities2}
\varrho_+=\begin{cases}R/\alpha,&\alpha>0,\\0,&\alpha=0,\end{cases}
 \qquad
 \varrho_-=\begin{cases}Q/(1-\alpha),&\alpha<1,\\0,&\alpha=1.\end{cases}
\end{equation}
In particular, $\varrho_{\pm}\left(  t,\cdot\right)  $ are measurable for each
time $t\in\left[  0,T\right]$, and 
\begin{equation}\label{eq:rho_pm}
R=\alpha\varrho_+,\qquad Q=(1-\alpha)\varrho_-.
\end{equation}

However, we cannot say anything about integrability of $\varrho_\pm$. We can only deduce some weighted estimate following from the uniform estimates \eqref{eq:energydefinitions}.  
To this end, for  $\gamma>1$, we consider the following convex extension of
$(a,x)\mapsto x^\gamma/a^{\gamma-1}$ defined on $a>0$ and $x\geq0$:
\begin{equation}
\tilde p_\gamma(x,a) =\left\{
\begin{array}{ll}
\frac{x^\gamma}{a^{\gamma-1}} &\text{ for } a>0,\ x\geq0,\\
0 &\text{ for } a=0,\ x=0,\\
+\infty &\text{ for } a=0,\ x>0.
\end{array} 
\right.
\end{equation}
Note that separation of phases again gives that 
\[\tilde p_{\gamma^+}(r_k,s_k)=r_k^\gamma=p_+(r_k).\]
On the other hand, due to definition \eqref{eq:phasedensities2}, we have
\[\tilde p_{\gamma^+}(R,\alpha)=\alpha p_+(\varrho_+).\]
Now, for every measurable $I\subset(0,T)$, application of
\cite[Theorem~11.27]{FeireislNovotny2017}, and uniform estimates \eqref{eq:energydefinitions} give
\begin{equation*}\label{eq:proposal_phase_power}
\int_I\!\int_\Omega\tilde p_{\gamma^+}(R,\alpha)\dxt
 \leq\liminf_k\int_I\!\int_\Omega \tilde p_{\gamma^+}(r_k,s_k)\dxt
 \leq C_T|I|.
\end{equation*}
Which implies \(\int_I\!\int_\Omega\alpha p_+(\varrho_+)\dxt\leq C_T|I|.\)
Since $I$ is arbitrary, we deduce that
\[\alpha p_+(\varrho_+)\in L^\infty(0,T;L^1(\Omega)),\]
The same argument applied to the minus phase implies
\[(1-\alpha) p_-(\varrho_-)\in L^\infty(0,T;L^1(\Omega)).\]

The same theorem applied to the closed convex function
$(a,x,y)\mapsto(x^2+y^2)/a$ and to
$(s_k,r_k,s_k\diver u_k)\rightharpoonup(\alpha,R,G)$ gives
\begin{equation}\label{eq:phaseWeighted_1}
 \int_{Q_T}\frac{R^2+G^2}{\alpha}\dxt
 \leq\liminf_k\int_{Q_T}
        [r_k^2+s_k|\diver u_k|^2]\dxt<\infty.
\end{equation}
Indeed, $s_k\in\{0,1\}$ and $r_k=s_kr_k$ imply
\[
 \frac{r_k^2+|s_k\diver u_k|^2}{s_k}
 =r_k^2+s_k|\diver u_k|^2,
\]
and the right-hand side is uniformly
bounded in $L^1(Q_T)$ by \eqref{eq:energybounds}--\eqref{eq:extra}.

\subsection{Renormalized equations}\label{sec:entropy}
We now recall the result on the renormalized continuity with the source term. Its proof, as a special case, can be found in  \cite[Theorem 11.36, Lemma 11.13]{FeireislNovotny2017}
\begin{lemma}[Renormalized transport equation]\label{lem:renormalization}
Suppose that the functions $(h,u,f)\in L^2_{loc}((0,T)\times\R^3)\times L^2(0,T;H^1_{loc}(\R^3, \R^3))\times L^2_{loc}((0,T)\times\R^3)$, where
$h\geq0$ a.e. in \((0,T)\times\R^3\) satisfy the transport equation
\begin{equation}\label{cont}
 \partial_t h+\diver(hu)=f\quad in\ \ {\cal D}'((0,T)\times \R^3),
\end{equation}
with the initial datum  \(h^0\in L^1_{loc}(\R^3)\), \(h^0\geq0\). Then
for any 
\begin{equation}
    b\in C^1([0,\infty])\cap W^{1,\infty}(0,\infty),
\end{equation}
we have
\begin{equation}\label{eq:sourcedrenorm}
 \partial_t b(h)+\diver(b(h)u)=fb'(h)+\left(b(h)-hb'(h)\right)\diver u,\quad in\ \ {\cal D}'((0,T)\times \R^3),
\end{equation}
with initial datum \(b(h^0)\). In particular, for every
\(\psi\in C_c^\infty([0,T)\times\overline\Omega)\),
\[
\begin{aligned}
&\int_{Q_T}
\Bigl[
 b(h)\partial_t\psi
 +b(h)u\cdot\nabla\psi
 +\bigl(b(h)-hb'(h)\bigr)\operatorname{div}u\,\psi
 +b'(h)f\,\psi
\Bigr]\,dx\,dt\\
&\qquad+\int_\Omega b(h^0)\psi(0)\,dx=0.
\end{aligned}
\]
If \(0\leq h,h^0\leq1\), the same conclusion holds for every \(b\in C^1([0,1])\).
\end{lemma}

Applying Lemma 3.1 with \(h=R\), \(u=U\), and \(f=0\), we obtain
\[
\partial_t b(R)+\operatorname{div}(b(R)U)
=\bigl(b(R)-Rb'(R)\bigr)\operatorname{div}U
\]
for every \(b\in C^1([0,\infty))\) with bounded derivative, with initial datum \(b(R^0)\). The same assertion holds for \(Q\).
Applying the lemma  with \(h=\alpha\), \(u=U\), and \(f=\overline{s\operatorname{div}u}\), gives
\[
\partial_t b(\alpha)+\operatorname{div}(b(\alpha)U)
=
b'(\alpha)\,\overline{s\operatorname{div}u}
+\bigl(b(\alpha)-\alpha b'(\alpha)\bigr)\operatorname{div}U,
\qquad b\in C^1([0,1]).
\]
This verifies the validity of the equations \eqref{eq:weakrenormalizedmass} and \eqref{eq:BNfractionweak}.

For later use, we also cite an extract of \cite[Proposition 3.1]{Novotny2020} giving joint renormalized equation for $(\alpha,R)$ with forcing.

\begin{proposition}[Joint renormalization]\label{prop:joint-renormalization}
Let $0\leq\alpha\leq1$, $R\geq0$, and suppose that
\[
 R\in L^2(Q_T)\cap L^\infty(0,T;L^{\gamma^+}(\Omega)),
 \quad
 U\in L^2(0,T;H^1_0(\Omega;\mathbb R^3)),
 \qquad G\in L^2(Q_T).
\]
Assume that
\[
 \partial_t\alpha+\diver(\alpha U)=G,
 \qquad
 \partial_tR+\diver(RU)=0
\]
 in $\mathcal D'(Q_T)$, with initial data $\alpha^0,R^0$.
Then, for every
\[
 B\in C^1([0,\infty)^2),
 \qquad
 \partial_aB,\partial_rB\in L^\infty([0,\infty)^2),
\]
we have
\begin{equation}\label{eq:joint-renormalization}
 \begin{aligned}
 \partial_t B(\alpha,R)
 +\diver\bigl(B(\alpha,R)U\bigr)
 ={}&\partial_aB(\alpha,R)\,G\\
 &+\bigl[B(\alpha,R)
       -\alpha\partial_aB(\alpha,R)
       -R\partial_rB(\alpha,R)\bigr]\diver U
 \end{aligned}
\end{equation}
in $\mathcal D'(Q_T)$.

Moreover,
\[
 B(\alpha,R)\in C([0,T];L^1(\Omega)),
 \qquad
 B(\alpha,R)|_{t=0}=B(\alpha^0,R^0).
\]
Identity \eqref{eq:joint-renormalization} also holds in the
time-integrated form up to the boundary.
\end{proposition}

\begin{proof}
Novotn\'y states the proposition for two comparable nonnegative quantities, that we construct below. 
Let $\xi=\alpha+R$ and $\zeta=\alpha+2R$. Then
$\xi\leq\zeta\leq2\xi$, and
\[
\partial_t\xi+\diver(\xi U)=G,
\qquad
 \partial_t\zeta+\diver(\zeta U)=G.
 \]
 We now apply \cite[Proposition~3.1(1), equation~(3.6)]{Novotny2020}
to $(\xi,\zeta)$, with $v=\diver U$ and
\[
 \widetilde B(\xi,\zeta)
 =B(2\xi-\zeta,\zeta-\xi).
\]
The comparison condition required there is satisfied with constants $1$ and $2$. Moreover,
\[
 \begin{aligned}
 \partial_\xi\widetilde B+\partial_\zeta\widetilde B
   &=\partial_aB(\alpha,R),\\
 \xi\partial_\xi\widetilde B+\zeta\partial_\zeta\widetilde B
   &=\alpha\partial_aB(\alpha,R)+R\partial_rB(\alpha,R),
 \end{aligned}
\]
where the derivatives of $\widetilde B$ are evaluated at
$(\xi,\zeta)=(\alpha+R,\alpha+2R)$.
\end{proof}


\subsection{The effective viscous-flux}
\label{sec:flux}
In this section we formulate result concerning the weak compactness of the effective viscous flux. Set
\begin{equation}\label{eq:effectivenotation}
 F_k=(2\mu+\lambda)\diver u_k-p_k,\qquad
 F=(2\mu+\lambda)\diver U-\overline p.
\end{equation}
The lemma below establishes essential fact that the product of $F_k$  with a bounded sequence satisfying transport equation with a $L^2$ source converges to the product of the weak limits in the sense of distributions.
This is a classical result for single-component Navier-Stokes equations, and its generalization to two-fluid algebraic pressure closure system was discussed, for example, in \cite[Proposition~12]{NovotnyPokorny2020}. We refer the reader to \cite[Proposition 7.36]{NovotnyStraskraba2004} where this result was proven under weaker, and somehow optimal, assumptions for the compressible non-steady Navier-Stokes equations.

\begin{lemma}[Effective viscous flux]\label{lem:eff_Flux}
Let $\beta=\min\{\gamma^+,\gamma^-\}$, let $p_0$ be defined by
\eqref{eq:compactexponents}, and let 
\(p_0^*=\frac{3p_0}{3-p_0}. \)
Let
\[
 1<q,r<\infty,\qquad
 1<v\leq\infty,\qquad
 \max\{2,r'\}<w\leq\infty,
\]
satisfy
\[
 \frac1{p_0}+\frac1q<\frac56,
 \qquad
 \frac1{p_0^*}+\frac1v<1,
\]
where $r'=r/(r-1)$.

Suppose that $(r_k,q_k,u_k)$ satisfy the continuity and momentum
equations \eqref{eq:weakmass} and \eqref{eq:weakmomentum}, and set
\[
 \varrho_k=r_k+q_k,\qquad
 \varrho=R+Q,\qquad
 p_k=p_+(r_k)+p_-(q_k).
\]
Assume that
\begin{equation}\label{subseq_ass}
\begin{aligned}
 \varrho_ku_k
 &\longrightarrow\varrho U
 &&\text{in }C_{\mathrm{weak}}
       ([0,T];L^{p_0}(\Omega;\mathbb R^3)),\\
 u_k
 &\rightharpoonup U
 &&\text{in }L^2(0,T;H^1_0(\Omega;\mathbb R^3)),\\
 r_k\longrightarrow R,\qquad
 q_k&\longrightarrow Q
 &&\text{in }C_{\mathrm{weak}}
       ([0,T];L^\beta(\Omega)),\\
 p_k
 &\rightharpoonup\overline p
 &&\text{in }L^r(Q_T),\\
 h_k
 &\longrightarrow h
 &&\text{in }C_{\mathrm{weak}}
       ([0,T];L^q(\Omega)),\\
 h_k
 &\rightharpoonup h
 &&\text{in }L^w(Q_T),\\
 f_k
 &\rightharpoonup f
 &&\text{in }L^v(Q_T).
\end{aligned}
\end{equation}
The last two convergences are understood weakly-* when
$w=\infty$ or $v=\infty$, respectively.
Suppose moreover that
\begin{equation}\label{eq:observablebalance}
 \partial_t h_k+\diver(h_ku_k)=f_k
 \qquad\text{in }\mathcal D'(Q_T).
\end{equation}

Then, after passing to a subsequence,
\begin{equation}\label{subseq}
\begin{aligned}
 h_k\diver u_k
 &\rightharpoonup\overline{h\diver u}
 &&\text{in }L^{\frac{2w}{2+w}}(Q_T),\\
 h_kp_k
 &\rightharpoonup\overline{hp}
 &&\text{in }L^{\frac{wr}{w+r}}(Q_T).
\end{aligned}
\end{equation}
Moreover,
\begin{equation}
 (2\mu+\lambda)\overline{h\diver u}-\overline{hp}
 =
 h\bigl((2\mu+\lambda)\diver U-\overline p\bigr)
 \qquad\text{a.e. in }Q_T.
\end{equation}
In particular, with $F_k,F$ defined in
\eqref{eq:effectivenotation},
\[
 h_kF_k\rightharpoonup hF
 \qquad\text{in }L^{p_{\min}}(Q_T),
 \qquad
 p_{\min}
 =\min\left\{\frac{2w}{2+w},\frac{wr}{w+r}\right\}>1.
\]
\end{lemma}

We will be using this lemma repetitively for functions $h_k$ involving combinations of $s_k$ and bounded bounded nonlinear functions of $r_k$ and $q_k$. This will eventually allow us to identify the limit $\overline p$. For additional restrictions on $\gamma^\pm$ we will also apply it to $h_k=r_k, q_k$.

\subsection{Application of the effective-viscus flux equality}

\begin{proposition}[Phase-weighted effective flux]\label{prop:phaseflux}
Let $p_{min}=\min\{{\frac{2w}{2+w}}, {\frac{wr}{w+r}} \}$, where $w,r$ are as in Lemma \ref{lem:eff_Flux}.
For every 
\[
 b\in C^1([0,\infty)),\qquad\text{with}\quad
 \|b\|_\infty+\|zb'(z)\|_\infty<\infty,
\]
we have
\begin{equation}\label{eq:mainflux}
 \begin{aligned}
 s_kb(\varrho_k)F_k
       &\rightharpoonup\overline{s b(\varrho)}\,F,\\
 (1-s_k)b(\varrho_k)F_k
       &\rightharpoonup\overline{(1-s)b(\varrho)}\,F
 \end{aligned}
 \quad\text{weakly in }L^{p_{min}}(Q_T).
\end{equation}

\end{proposition}
\begin{proof}
Nothe that, since $s_k\in\{0,1\}$ and $r_k=s_k\varrho_k$, we have
\begin{equation}\label{eq:phaseidentity}
 s_kb(\varrho_k)=b(r_k)-b(0)(1-s_k).
\end{equation}
Indeed, both sides equal $b(\varrho_k)$ when $s_k=1$ and zero when $s_k=0$.

Subtracting $b(0)$ times the equation for $1-s_k$ from the renormalized
$r_k$ equation \eqref{eq:weakrenormalizedmass}, and using
$r_kb'(r_k)=s_k\varrho_kb'(\varrho_k)$, we obtain
\begin{equation}\label{eq:phasebalance}
 \partial_t(s_kb(\varrho_k))
 +\diver(s_kb(\varrho_k)u_k)
 =s_k[b(\varrho_k)-\varrho_kb'(\varrho_k)]\diver u_k,
\end{equation}
in the weak sense. 

Next, take $h_k=s_kb(\varrho_k)$ and let $f_k$ be the right-hand side of
\eqref{eq:phasebalance}. Choose in Lemma~\ref{lem:eff_Flux}
\[
 r=1+\delta,\qquad w=\infty,\qquad v=2,
 \qquad q>\frac{6\beta}{2\beta-3}.
\]
Then $p_{\min}=\min\{2,1+\delta\}>1$. Due to assumptions on $b$, we have 
\[
 \|h_k\|_\infty\leq\|b\|_\infty,\qquad
 \|f_k\|_{L^2}\leq
       \|b-zb'\|_\infty\|\diver u_k\|_{L^2}.
\]
Moreover, using \eqref{eq:phasebalance} we can estimate $\partial_t h_k$ in
$L^2(0,T;H^{-1}(\Omega))$. The compactness argument used
for \eqref{conv_faible} therefore gives
\[
 h_k\longrightarrow\overline{s b(\varrho)}
 \quad\text{in }C_{\rm weak}([0,T];L^q(\Omega)).
\]
Thus all the assumptions of Lemma~\ref{lem:eff_Flux} are satisfied, and the first convergence in \eqref{eq:mainflux} follows. For the second one observe that $(1-s_k)b(\varrho_k)=b(q_k)-b(0)s_k$,
and proceed similarly as before.
\end{proof}

\section{Identification of pressure}
\label{sec:phase}

\subsection{Heuristic explanation}
To explain  the mechanism of the proof, we first present our argument for the compactness of the pressure schematically.
It can be made rigorous under additional restriction on the pressure exponents:
\begin{equation}\label{eq:strongergamma}
 \beta=\min\{\gamma^+,\gamma^-\}>\frac 92.
\end{equation}
Under this constraint, the Bogovskii type argument yields the stronger estimate
\begin{equation}\label{eq:strongergain}
 \sup_k\int_{Q_T}
 \bigl(r_k^{(1+\sigma)\gamma^+}
       +q_k^{(1+\sigma)\gamma^-}\bigr)\dxt\leq C_T\quad  \text{with}\quad \frac1\beta<\sigma<\frac13,
\end{equation}
which implies the following uniform bounds
\begin{equation}\label{eq:productintegrability}
 \begin{aligned}
 p_k&\text{ bounded in }L^{1+\sigma}(Q_T),\qquad
 \varrho_k\text{ bounded in }L^{\beta(1+\sigma)}(Q_T),\\
 \varrho_kp_k&\text{ bounded in }L^a(Q_T),\qquad
 a:=\frac{\beta(1+\sigma)}{\beta+1}>1.
 \end{aligned}
\end{equation}
We first note that due to the uniform estimate of $r_k$ in $L^2(Q_T)$ DiPerna-Lions transport theory applies to the continuity equation for $r_k$. It implies the following renormalized equation
\[
 \partial_t(r_k\log r_k)+\diver(r_k\log r_k\,u_k)=-r_k\diver u_k\quad \text{in } \mathcal D'((0,T)\times\R^3).
\]
The compactness argument of \eqref{eq:convectiveclosure} allows us then to pass to the weak limit to justify
\begin{equation}\label{eq:meanentropy}
 \partial_t\overline{r\log r}+\diver(\overline{r\log r}\,U)
 =-\overline{r\diver u} \quad \text{in } \mathcal D'((0,T)\times\R^3).
\end{equation}
We now compare it with the joint renormalized equation for the limiting quantities $\alpha$ and $R$. Using the definition of $\varrho_+=\frac{R}{\alpha}$, we eventually 
\begin{equation}\label{eq:perspectivebalance}
 \partial_t(\alpha\varrho_+\log\varrho_+)
 +\diver(\alpha\varrho_+\log\varrho_+\,U)=-\varrho_+G\quad \text{in } \mathcal D'((0,T)\times\R^3),
\end{equation}
which holds with the same initial data as in \eqref{eq:meanentropy}. Recall that we dented $G=\overline{s\diver u}$.
We now consider 
$$\mathcal J:=\overline{r\log r}-\alpha\varrho_+\log\varrho_+,$$
that will measure the evolution of  $\alpha$--weighted  oscillations for the density. Subtracting \eqref{eq:perspectivebalance}
from \eqref{eq:meanentropy}, we obtain the equation $\mathcal J$ satisfies
\begin{equation}\label{eq:J}
 \partial_t\mathcal J+\diver(\mathcal J U)
 =-\overline{r\diver u}+\varrho_+\overline{s\diver u},\qquad \mathcal J|_{t=0}=0.
\end{equation}

We now use the effective viscous flux identity to show that the r.h.s. of equation \eqref{eq:J} has a sign. First, taking $b=1$ in
\eqref{eq:mainflux} we get $\overline{s F}=\alpha F$. Recalling notation \eqref{eq:effectivenotation} and that separation of phases $r_k q_k=0$ implies $s_kp_k=p_+(r_k)$, we get
\begin{equation}\label{eq:baseflux}
 (2\mu+\lambda)\overline{s \diver u}-\overline{p_+(r)}=\alpha F.
\end{equation}

Now we take $b=Id$ in \eqref{eq:mainflux} (which is not admissible, and thus the argument is formal). Since separation of phases implying $
r_kp_k=r_kp_+(r_k)+r_kp_-(q_k)= r_kp_+(r_k)$, we obtain
\begin{equation}\label{eq:densityflux}
 \begin{aligned}
 (2\mu+\lambda)\overline{r\diver u}-\overline{r p_+(r)}
       &=RF.
 \end{aligned}
\end{equation}

Substituting $\varrho_+\times$\eqref{eq:baseflux} -\eqref{eq:densityflux} into the
equation \eqref{eq:J}, the terms involving $F$ cancel since
$R=\alpha\varrho_+$. Hence
\begin{equation}\label{eq:J_heur}
 \partial_t\mathcal J+\diver(\mathcal J U)
 =-\frac{\overline{r p_+(r)}-\varrho_+\overline{p_+(r)}}{2\mu+\lambda},\qquad \mathcal J|_{t=0}=0.
\end{equation}
Formally, integrating this equality over $Q_T$, using monotonicity of $s\mapsto s^{\gamma_+}$ and convexity of $x\mapsto x\log x$, one
obtains
\begin{equation}\label{eq:entropyidentified}
 \overline{r\log r}=\alpha\varrho_+\log\varrho_+
 \qquad\text{a.e. in }Q_T.
\end{equation}
We now want to turn the r.h.s. to the convex function of one phase. This could be achieved by adding $(1-\alpha)\varrho_+\log\varrho_+$, which could be obtained as a limit of $(1-s_k)\varrho_+\log\varrho_+$. Observe now that we have
\begin{equation}\label{loglog}
r_k\log r_k+(1-s_k)\varrho_+\log\varrho_+=z_k\log z_k,
\end{equation}
where $z_k$ is given\footnote{The construction of $z_k$ is motivated by the following informal consideration: if we assume the following two-scale development ansatz of the unknowns $s_k(t,x)=\chi(t,x,kx)+o(1),\text{ } \varrho_k(t,x)=\chi(t,x,kx)\varrho_+(t,x)+(1-\chi(t,x,kx))\varrho_-(t,x)+ o(1)$ as $k\rightarrow \infty$ then we find that $r_k(t,x)+(1-s_k(t,x))\varrho_+(t,x)= \varrho_+(t,x) + o(1)$ as $k\rightarrow \infty$.} as $k\rightarrow\infty$ by
\[
z_k=r_k+(1-s_k)\varrho_+
\]
Indeed, it follows from the separation of phases again. If \(s_k=1\), then \(z_k=r_k\) and \(1-s_k=0\), thus \(r_k\log r_k=z_k\log z_k\).  If \(s_k=0\), segregation gives \(r_k=0\), so \(z_k=\varrho_+\), thus \(\varrho_+\log\varrho_+=z_k\log z_k.\)
Passing to the weak limit for $k\to\infty$ can the justify that
\[
z_k\rightharpoonup R+(1-\alpha)\varrho_+
=\alpha\varrho_++(1-\alpha)\varrho_+
=\varrho_+,
\]
while passing to the limit $k\to\infty$ in \eqref{loglog}, one obtains
\[
\begin{aligned}
z_k\log z_k
&\rightharpoonup
\overline{r\log r}+(1-\alpha)\varrho_+\log\varrho_+
=\varrho_+\log\varrho_+.
\end{aligned}
\]
Using convexity of $x\log x$ again, criterion of \cite[Theorem~11.27]{FeireislNovotny2017}, eventually yields
\begin{equation}\label{zrho0}
z_k\longrightarrow\varrho_+
\quad\text{a.e. on }Q_T.
\end{equation}
This point-wise convergence can be then used to identify the limits of all nonlinear functions of the density, involving the pressure and the r.h.s. of the equation for $\alpha$.



\subsection{Rigorous argument for  \texorpdfstring{$\gamma^\pm\geq9/5$}{of general exponents}}
We now explain how the above argument can be made rigorous in the full range of admissible pressure exponents, without the restriction \eqref{eq:strongergamma}. Similar constructions can be found in \cite{FeireislNovotny2017}. We will prove the following result.
\begin{proposition}[Identification within each phase]\label{prop:phasemarginal}
 For every $b\in C_b([0,\infty))$,
\begin{equation}\label{eq:mainatoms}
 \overline{s b(\varrho)}=\alpha b(\varrho_+),\qquad
 \overline{(1-s)b(\varrho)}=(1-\alpha)b(\varrho_-).
\end{equation}
Here the left-hand sides denote weak-* limits in $L^\infty(Q_T)$.
\end{proposition}

\begin{proof}
We prove the first identity; the other phase is treated in the same way.
\paragraph{Nonnegativity of truncated measure of oscillations $\mathcal J_n$.} 
Recall that $G=\overline{s\diver u}$ and that
$R=G=0$ on $\{\alpha=0\}$ by \eqref{eq:phaseWeighted_1}.
Let
\[
 T_n(z)=\frac{nz}{n+z},\qquad
 L_n(z)=z\log z-z\log(1+z/n),\qquad L_n(0)=0.
\]
The function $T_n$ is bounded, increasing and concave, while $L_n$
is strictly convex and has at most linear growth. Moreover,
\begin{equation}\label{LL'}
 zL_n'(z)-L_n(z)=T_n(z),\qquad
 L_n''(z)=\frac{n^2}{z(n+z)^2}>0\quad(z>0).
\end{equation}
Segregation and the values at zero give
\begin{equation}\label{segregate}
 \begin{gathered}
 L_n(r_k)=s_kL_n(\varrho_k),\qquad
 T_n(r_k)=s_kT_n(\varrho_k),\\
 T_n(r_k)p_k=T_n(r_k)p_+(r_k),\qquad
 s_kp_k=p_+(r_k).
 \end{gathered}
\end{equation}
For fixed $c>0$, convexity of $L_n(\varrho_k)$ gives
\[L_n(\varrho_k)\geq L_n(c)+L'(c)(\varrho_k-c)=L'(c)\varrho_k-T_n(c),\]
where we have used \eqref{LL'}. Multiplying both sides by $s_k\geq0$ and using \eqref{segregate} gives
\begin{equation}\label{before_limit}
    L_n(r_k)\geq L_n'(c)r_k-T_n(c)s_k.
\end{equation}
The function $L_n(r_k)$ has at most linear growth as
\begin{equation}\label{linear}
|L_n(r_k)|\le n+|\log n|\,r_k,
\end{equation}
thus for fixed $n$ it is bounded in $L^2(Q_T)$ uniformly w.r.t. $k$. 

For a fixed $c\in\mathbb{Q},\text{ }c>0$ and every nonnegative
$\varphi\in C_c^\infty(Q_T)$, the weak convergences imply
\[
 \begin{aligned}
 &\int_{Q_T}\varphi
 \bigl[\overline{L_n(r)}-L_n'(c)R+T_n(c)\alpha\bigr]\dxt\\
 &\qquad=
 \lim_{k\to\infty}\int_{Q_T}\varphi
 \bigl[L_n(r_k)-L_n'(c)r_k+T_n(c)s_k\bigr]\dxt
 \geq0.
 \end{aligned}
\]
Consequently,
\begin{equation}
\label{eq:pozitiv}
 \overline{L_n(r)}
 \geq L_n'(c)R-T_n(c)\alpha
 \qquad\text{on } Q_T-S_c,
\end{equation}
where $S_c$ has Lebesgue measure $0$. Using that the rationals are countable we conclude that there exists a measurable set $S$ with Lebesgue measure $0$ such that \eqref{eq:pozitiv} holds true for all rationals $c$ and all $(t,x)\in Q_T-S$. Using a density argument, we extend \eqref{eq:pozitiv} a.e. on $Q_T$ for any $c\geq0$.
and the continuity of the r.h.s. as a function of $c\in(0,,\infty)$ allow us to take all $c>0$. We may now choose $c$ depending on $(t,x)$.
\begin{itemize}
    \item On $\{\alpha>0,\ R>0\}$, take $c=R/\alpha=\varrho_+$. Using
$zL_n'(z)-L_n(z)=T_n(z)$, we obtain
\[
 \begin{aligned}
 \overline{L_n(r)}
 &\geq L_n'(\varrho_+)R-\alpha T_n(\varrho_+)\\
 &=\alpha\bigl[\varrho_+L_n'(\varrho_+)-T_n(\varrho_+)\bigr]
 =\alpha L_n(\varrho_+).
 \end{aligned}
\]
\item On $\{\alpha>0,\ R=0\}$, we have $\varrho_+=0$, and the inequality
reduces to
\[
 \overline{L_n(r)}\geq-\alpha T_n(c).
\]
Letting $c\to0$ gives
$\overline{L_n(r)}\geq0=\alpha L_n(\varrho_+)$.
\item Finally, on $\{\alpha=0\}$ we already know that $R=0$.
Thus the limiting inequality gives
$\overline{L_n(r)}\geq0$ for any fixed $c>0$.
With the convention that $\alpha L_n(\varrho_+)=0$ on this set,
we conclude that
\begin{equation}\label{eq:entropydefect}
 \mathcal J_n:=\overline{L_n(r)}
              -\alpha L_n(\varrho_+)\geq0
 \qquad\text{a.e. in }Q_T.
\end{equation}
\end{itemize}

\paragraph{The equation for $\mathcal J_n$.}
We now derive the equation for $\mathcal J_n$. The idea is very similar as before: we first renormalize equation for $r_k$ and pass to the limit $k\to\infty$ and compare it with joint renormalized equation for $\alpha$ and $R$. The latter is obtained thanks to Proposition \ref{prop:joint-renormalization}.

For a fixed $\varepsilon>0$, we apply \eqref{eq:weakrenormalizedmass} with the function
$b_\varepsilon(z)=L_n(z+\varepsilon)-L_n(\varepsilon)\in C^1([0,\infty)) \text{ and } b_\varepsilon\in L^{\infty}([0,\infty)])$. For every
\(
\psi\in C_c^\infty([0,T)\times\overline\Omega)\),
we therefore have
\begin{equation}\label{b_ep}
    \begin{aligned}
&\int_{Q_T}
 b_\varepsilon(r_k)
 \bigl(\partial_t\psi+u_k\cdot\nabla\psi\bigr)
 \,\mathrm dx\,\mathrm dt+
 \int_{Q_T}
 \bigl[b_\varepsilon(r_k)-r_kb_\varepsilon'(r_k)\bigr]
 \operatorname{div}u_k\,\psi
 \,\mathrm dx\,\mathrm dt\\
&\quad+
 \int_\Omega b_\varepsilon(r_k^0)\psi(0)\,\mathrm dx
 =0.
\end{aligned}
\end{equation}
We now want to let $\varepsilon\to 0$ in this equation.
First note that, by continuity of \(L_n\) on \([0,\infty)\), we have \(
b_\varepsilon(r_k)\longrightarrow L_n(r_k)\) as $\varepsilon\to 0$, and so it is uniformly bounded in $L^2(Q_T)$.
The linear-growth estimate \eqref{linear}
also yields
\[
|b_\varepsilon(r_k)|
\le 2n+|\log n|(r_k+2)
\le C_n(1+r_k),
\]
which gives an integrable majorant for
\[
\left|
b_\varepsilon(r_k)
(\partial_t\psi+u_k\cdot\nabla\psi)
\right|
\le C_{n,\psi}(1+r_k)(1+|u_k|)\in L^1(Q_T).
\]
This allows to pass to the limit in the first integral of \eqref{b_ep}.

For the second integral, we first check that
\begin{equation}\label{dominated}
    \begin{aligned}
r_kb_\varepsilon'(r_k)-b_\varepsilon(r_k)
&=\int_0^{r_k} t\,b_\varepsilon''(t)\,\mathrm dt=\int_0^{r_k}
 \frac{n^2t}{(t+\varepsilon)(n+t+\varepsilon)^2}\,\mathrm dt,
\end{aligned}
\end{equation}
where we used \(b_\varepsilon(0)=0\).
Consequently, we can estimate
\[
0\le r_kb_\varepsilon'(r_k)-b_\varepsilon(r_k)
\le \int_0^{r_k}\frac{n^2}{(n+t)^2}\,\mathrm dt\le n.
\]
Dominated convergence theorem applied to \eqref{dominated} together with \(\int_0^{r_k}\frac{n^2}{(n+t)^2}\,\mathrm dt
=T_n(r_k)\) therefore implies
\[
b_\varepsilon(r_k)-r_kb_\varepsilon'(r_k)
\longrightarrow -T_n(r_k)
\qquad\text{for every fixed }r_k\ge0.
\]
Because
\[
\left|
[b_\varepsilon(r_k)-r_kb_\varepsilon'(r_k)]
\operatorname{div}u_k\,\psi
\right|
\le n\|\psi\|_\infty|\operatorname{div}u_k|,
\]
is integrable, one can pass to the limit $\varepsilon\to 0$ also in the second integral of \eqref{b_ep}.
Finally, the initial integrand satisfies
\[
|b_\varepsilon(r_k^0)\psi(0)|
\le C_n\|\psi(0)\|_\infty(1+r_k^0)\in L^1(\Omega),
\]
and converges pointwise to \(L_n(r_k^0)\psi(0)\).
Ultimately, letting  $\varepsilon\to 0$ in \eqref{b_ep}, and rearranging the terms results in:
\[
\begin{aligned}
&\int_{Q_T}
L_n(r_k)\bigl(\partial_t\psi+u_k\cdot\nabla\psi\bigr)
\,\mathrm dx\,\mathrm dt
+\int_\Omega L_n(r_k^0)\psi(0)\,\mathrm dx=
\int_{Q_T}T_n(r_k)\operatorname{div}u_k\,\psi
\,\mathrm dx\,\mathrm dt,
\end{aligned}
\]
for every \(\psi\in C_c^\infty([0,T)\times\overline\Omega)\).

We now want to let $k\to\infty$ in the above.
For fixed $n$, $L_n(r_k)$ is bounded in
$L^\infty(0,T;L^{\gamma^+}(\Omega))$, its transport flux is bounded in
$L^2(0,T;L^{6\gamma^+/(\gamma^++6)}(\Omega))$, and its source is bounded
in $L^2(Q_T)$. The compactness argument of
\eqref{eq:convectiveclosure}, using $\gamma^+>6/5$, gives strong
convergence of $L_n(r_k)$ in $C([0,T];H^{-1}(\Omega))$ and identifies
the limiting flux as $\overline{L_n(r)}U$.
Since
\[
 L_n(r_k^0)=\chi_kL_n(\varrho_+^0),
\]
the initial term converges to $\alpha^0L_n(\varrho_+^0)$, by
\eqref{eq:partitionlimit} and
$L_n(\varrho_+^0)\in L^1(\Omega)$.
Letting $k\to\infty$, we obtain
\begin{equation}\label{Jn1}
    \begin{aligned}
&\int_{Q_T}\overline{L_n(r)}
   \bigl(\partial_t\psi+U\cdot\nabla\psi\bigr)\,\mathrm dx\,\mathrm dt
 +\int_\Omega\alpha^0L_n(\varrho_+^0)\psi(0)\,\mathrm dx=
 \int_{Q_T}\overline{T_n(r)\operatorname{div}u}\,
       \psi\,\mathrm dx\,\mathrm dt,
\end{aligned}
\end{equation}
for every \(\psi\in C_c^\infty([0,T)\times\overline\Omega)\).

We now want to derive the equation for $\alpha L_n(\varrho_+)$. 
To this purpose we apply Proposition~\ref{prop:joint-renormalization} with the pair  $(\alpha,R)=(\alpha,\alpha\varrho_+)$, with the approximate function
\[
B_\varepsilon(\alpha,\alpha\varrho_+)
=
(\alpha+\varepsilon)
L_n\!\left(\frac{\alpha\varrho_++\varepsilon}
                 {\alpha+\varepsilon}\right)
-\varepsilon L_n(1).
\]
The identity \(zL_n'(z)-L_n(z)=T_n(z)\) gives
\[
\begin{aligned}
\partial_1B_\varepsilon(\alpha,\alpha\varrho_+)
&=-T_n\!\left(
       \frac{\alpha\varrho_++\varepsilon}{\alpha+\varepsilon}
       \right),\qquad
\partial_2B_\varepsilon(\alpha,\alpha\varrho_+)
=L_n'\!\left(
       \frac{\alpha\varrho_++\varepsilon}{\alpha+\varepsilon}
       \right).
\end{aligned}
\]
Since \(0\le\alpha\le1\) and \(\alpha\varrho_+\ge0\), we have 
\(
\frac{\alpha\varrho_++\varepsilon}{\alpha+\varepsilon}
\ge\frac{\varepsilon}{1+\varepsilon}
\ge\frac{\varepsilon}{2}
\). This, together with \(0\le T_n\le n\) and the formula for \(L_n'\), this implies
\[
\begin{aligned}
\left|\partial_1B_\varepsilon(\alpha,\alpha\varrho_+)\right|
&\le n,\qquad 
\left|\partial_2B_\varepsilon(\alpha,\alpha\varrho_+)\right|
\le C_n(1+|\log\varepsilon|).
\end{aligned}
\]
Thus for fixed $n$ and $\varepsilon$ the derivatives are bounded  and $B_\varepsilon$ is an admissible renormalization, and Proposition~\ref{prop:joint-renormalization} gives
\begin{equation}\label{weakBe}
\begin{aligned}
&\int_{Q_T}
B_\varepsilon(\alpha,\alpha\varrho_+)
\bigl(\partial_t\psi+U\cdot\nabla\psi\bigr)
\,\mathrm dx\,\mathrm dt
+\int_\Omega
B_\varepsilon(\alpha^0,\alpha^0\varrho_+^0)\psi(0)\,\mathrm dx
\\
&\quad=
\int_{Q_T}
T_n\!\left(
\frac{\alpha\varrho_++\varepsilon}{\alpha+\varepsilon}
\right)G\psi\,\mathrm dx\,\mathrm dt
\\
&\qquad-
\varepsilon\int_{Q_T}
\left[
L_n'\!\left(
\frac{\alpha\varrho_++\varepsilon}{\alpha+\varepsilon}
\right)
-
T_n\!\left(
\frac{\alpha\varrho_++\varepsilon}{\alpha+\varepsilon}
\right)
-L_n(1)
\right]
\operatorname{div}U\,\psi\,\mathrm dx\,\mathrm dt,
\end{aligned}
\end{equation}
for every \(\psi\in C_c^\infty([0,T)\times\overline\Omega)\).

We now let $\varepsilon\to 0$ in the formula above.
First, by \eqref{linear}, we may bound
\begin{equation}\label{estB}
\begin{aligned}
&\left|B_\varepsilon(\alpha,\alpha\varrho_+)\right|=\left|
(\alpha+\varepsilon)
L_n\!\left(
\frac{\alpha\varrho_++\varepsilon}{\alpha+\varepsilon}
\right)-\varepsilon L_n(1)
\right|\\
&\qquad\le
n(\alpha+\varepsilon)
+|\log n|(\alpha\varrho_++\varepsilon)
+\varepsilon|L_n(1)|
\le C_n(1+\alpha\varrho_+).
\end{aligned}
\end{equation}
On \(\{\alpha>0\}\), the expression inside the absolute value converges as $\varepsilon\to0$ to
\(\alpha L_n(\varrho_+)\). On \(\{\alpha=0\}\), it is exactly
\(
\varepsilon L_n(1)-\varepsilon L_n(1)=0.\)
Since \(\alpha\varrho_+,U\in L^2(Q_T)\), the first integral in \eqref{weakBe} is bounded, and one can use the dominate convergence theorem to let $\varepsilon\to0$ in this term.
For the first term on the right-hand side of \eqref{weakBe} we hae
\[
T_n\!\left(
\frac{\alpha\varrho_++\varepsilon}{\alpha+\varepsilon}
\right)G
\longrightarrow T_n(\varrho_+)G
\qquad\text{a.e. in }Q_T.
\]
Indeed, on \(\{\alpha>0\}\) this follows from continuity of \(T_n\). On \(\{\alpha=0\}\), both products vanish because \(G=0\), see \eqref{eq:phaseWeighted_1}. Moreover,
\[
\left|
T_n\!\left(
\frac{\alpha\varrho_++\varepsilon}{\alpha+\varepsilon}
\right)G
\right|
\le n|G|\in L^2(Q_T),
\]
so the dominated convergence theorem can be used to pass to the limit in this term as well.
The remaining integral on the r.h.s. of\eqref{weakBe} tends to zero. Indeed, the estimate  of $L_n'=\partial_2B$ from above together with the bound \(T_n\le n\) allow to estimate 
\begin{align*}
& \varepsilon \int_{Q_{T}}\left\vert L_{n}^{\prime }\!\left( \frac{\alpha
\varrho _{+}+\varepsilon }{\alpha +\varepsilon }\right) -T_{n}\!\left( \frac{%
\alpha \varrho _{+}+\varepsilon }{\alpha +\varepsilon }\right)
-L_{n}(1)\right\vert \diver U\,\psi \,\mathrm{d}x\,\mathrm{d}t\newline
\\
& \leq C_{n}\varepsilon (1+|\log \varepsilon |)\Vert \psi \Vert _{L^{\infty
}(Q_{T})}\Vert \diver U\Vert _{L^{1}(Q_{T})}\longrightarrow 0.
\end{align*}
Finally, $B_\varepsilon(\alpha^0,\alpha^0\varrho_+^0)$ converges pointwise and is integrable. Thus, letting $\varepsilon\to0 $ in \eqref{weakBe} we eventually obtain
\begin{equation}\label{Jn2}
    \begin{aligned}
&\int_{Q_T}\alpha L_n(\varrho_+)
   \bigl(\partial_t\psi+U\cdot\nabla\psi\bigr)\,\mathrm dx\,\mathrm dt
 +\int_\Omega\alpha^0L_n(\varrho_+^0)\psi(0)\,\mathrm dx=
 \int_{Q_T}T_n(\varrho_+)G\psi\,\mathrm dx\,\mathrm dt,
\end{aligned}
\end{equation}
for every \(\psi\in C_c^\infty([0,T)\times\overline\Omega)\).

Subtracting \eqref{Jn2} from \eqref{Jn1}, we therefore obtain
\begin{equation}\label{eq:J_n}
\begin{aligned}
&\int_{Q_T}\mathcal J_n
       \bigl(\partial_t\psi+U\cdot\nabla\psi\bigr)\,
       \mathrm dx\,\mathrm dt=
 \int_{Q_T}
 \left[\overline{T_n(r)\operatorname{div}u}
       -T_n(\varrho_+)\overline{s\diver u}\right]
       \psi\,\mathrm dx\,\mathrm dt.
\end{aligned}
\end{equation}
\paragraph{The sign of the r.h.s.}

To identify the right-hand side of \eqref{eq:J_n}, we use the effective viscous flux identity. Note that since $T_n$ and $T_n-zT_n'=nz^2/(n+z)^2$ are bounded, we can apply \eqref{eq:mainflux} with  $T_n$, and then with $1$. This gives
\[
 \begin{aligned}
 (2\mu+\lambda)\overline{T_n(r)\diver u}
       -\overline{T_n(r)p_+(r)}
       &=\overline{T_n(r)}\,F,\\
 (2\mu+\lambda)\overline{s\diver u}-\overline{p_+(r)}&=\alpha F.
 \end{aligned}
\]
\indent On $\{\alpha=0\}$, according to \eqref{eq:compatibility_func_of_r_1}-\eqref{eq:compatibility_func_of_r_2} and to the previous identities
\[
\overline{T_n(r)\diver u} = \overline{s\diver u}=0.
\]

On $\{\alpha>0\}$, multiply the second identity by
$\overline{T_n(r)}/\alpha$ and subtract it from the first. The terms
containing $F$ cancel, and we obtain
\[
 (2\mu+\lambda)\left(
   \overline{T_n(r)\diver u}
   -\frac{\overline{T_n(r)}}{\alpha}\overline{s\diver u}\right)
 =\overline{T_n(r)p_+(r)}
   -\frac{\overline{T_n(r)}}{\alpha}\overline{p_+(r)},
\]
where we used \eqref{segregate}.
Substituting it into equation  \eqref{eq:J_n} yields
\begin{equation}\label{eq:truncateddefect}
\begin{aligned}
 \partial_t \mathcal J_n+\diver(\mathcal J_nU)
 ={}&-\frac{\boldsymbol{1}_{\{\alpha>0\}}}{2\mu+\lambda}\left(
    \overline{T_n(r)p_+(r)}
         -\frac{\overline{T_n(r)}}{\alpha}\,\overline{p_+(r)}\right)\\
 &+\boldsymbol{1}_{\{\alpha>0\}}\left(T_n(\varrho_+)-\frac{\overline{T_n(r)}}{\alpha}\right)\overline{s\diver u}.
\end{aligned}
\end{equation}

Note that the first parenthesis in \eqref{eq:truncateddefect} is nonnegative. Indeed, passing to the weak limits in
\[
 s_k\bigl(T_n(\varrho_k)-T_n(c)\bigr)
       \bigl(p_+(\varrho_k)-p_+(c)\bigr)\geq0
\]
for some fixed $c\in\mathbb{Q},\text{ } c\geq0$ gives
\begin{equation}\label{eq:Cebasev}
 \overline{T_n(r)p_+(r)}
 -p_+(c)\overline{T_n(r)}
 -T_n(c)\overline{p_+(r)}
 +\alpha T_n(c)p_+(c)\geq0,
\end{equation}
for all $(t,x)\in Q_T-S_c $ where $S_c$ has Lebesgue measure $0$. Using that the rationals are countable we conclude that there exists a measurable set $S$ with Lebesgue measure $0$ such that \eqref{eq:Cebasev} holds true for all rationals $c$ and all $(t,x)\in Q_T-S$. Using a density argument, we extend \eqref{eq:Cebasev} a.e. on $Q_T$ for any $c\in \mathbb{R}$, $c\geq0$.

Then, at each point where $\alpha>0$, we choose
$c=\left(\frac{\overline{p_+(r)}}{\alpha}\right)^{\frac1{\gamma^+}}$ to cancel the last two terms in \eqref{eq:Cebasev}, and the first two give
\[
 \overline{T_n(r)p_+(r)}
 -\frac{\overline{T_n(r)}}{\alpha}\overline{p_+(r)}\geq0.
\]
Thus \eqref{eq:truncateddefect} now is replaced by the inequality
\begin{equation}\label{eq:entropyineq}
 \partial_t\mathcal J_n+\diver(\mathcal J_nU)
 \leq \boldsymbol{1}_{\{\alpha>0\}} \left[T_n(\varrho_+)-\frac{\overline{T_n(r)}}{\alpha}\right]
         \overline{s\diver u}.
\end{equation}
\paragraph{Vanishing of the oscillation defect.}
It remains to estimate the r.h.s. of \eqref{eq:entropyineq}.
Since
\[
 0\leq r_k-T_n(r_k)=\frac{r_k^2}{n+z}\leq\frac{r_k^2}{n}\in L^1(Q_T),
\]
thanks to \eqref{eq:extra}. Passing to the limit $k\to\infty$, restricting perhaps to a subsequence, we obtain 
\eqref{eq:extra} gives
\[
 0\leq\int_{Q_T}\bigl(R-\overline{T_n(r)}\bigr)\dxt
     \leq C_T/n.
\]
Since $T_n(r_k)$ increases with $n$, the weak limits
$\overline{T_n(r)}$ also increase along the common subsequence chosen
above. The last estimate therefore gives
$\overline{T_n(r)}\to R$ almost everywhere as $n\to\infty$, by a monotone convergence theorem. 
Let us now observe that the r.h.s. of \eqref{eq:entropyineq} is also integrable. Indeed, on $\{\alpha>0\}$, both
$T_n(\varrho_+)$ and $\overline{T_n(r)}/\alpha$ belong to
$[0,\varrho_+]$ and converge to $\varrho_+$. On the other hand, we have
\[
 \int_{Q_T}\varrho_+|\overline{s\diver u}|\dxt= \int_{Q_T}\varrho_+|G|\dxt
 \leq\left(\int_{Q_T}\frac{R^2}{\alpha}\dxt\right)^{1/2}
      \left(\int_{Q_T}\frac{G^2}{\alpha}\dxt\right)^{1/2}<\infty,
\]
by \eqref{eq:phaseWeighted_1}. Integrating \eqref{eq:entropyineq}, with zero
initial defect, and using dominated convergence thus gives
\[
 0\leq\mathop{\rm ess\,sup}_{0<\tau<T}\int_\Omega \mathcal J_n(\tau)\dx
 \leq\int_{Q_T}\left|T_n(\varrho_+)
        -\frac{\overline{T_n(r)}}{\alpha}\right||G|\dxt
 \longrightarrow0.
\]

Finally, let us also observe that
\[
 0\leq r_k\log r_k-L_n(r_k)
 =r_k\log(1+r_k/n)\leq\frac{r_k^2}{n}\in L^1(Q_T).
\]
 The $L^2(Q_T)$ bound on
$r_k$ gives uniform integrability of $r_k\log r_k$, hence passing to the weak limits we obtain
\[
 \bigl\|\overline{r\log r}-\overline{L_n(r)}\bigr\|_{L^1(Q_T)}
 \leq\frac1n\sup_k\|r_k\|_{L^2(Q_T)}^2\longrightarrow0.
\]
Also $\alpha\varrho_+^2\in L^1(Q_T)$ by \eqref{eq:phaseWeighted_1}, so
\[
 \bigl\|\alpha\varrho_+\log\varrho_+-\alpha L_n(\varrho_+)\bigr\|_{L^1(Q_T)}
 \leq\frac1n\int_{Q_T}\alpha\varrho_+^2\dxt\longrightarrow0.
\]
Together with
$\mathcal J_n\to0$ in $L^1(Q_T)$, this proves
\eqref{eq:entropyidentified}, namely
\[
 \overline{r\log r}=\alpha\varrho_+\log\varrho_+
       \quad\text{a.e. in }Q_T.
\]

\paragraph{Identification of the nonlinearities.}
We use the same functions as in the heuristic argument, namely
\[
z_k=r_k+(1-s_k)\varrho_+,
\]
which due to segregation of phases yields
\begin{equation}\label{loglog_2}
r_k\log r_k+(1-s_k)\varrho_+\log\varrho_+=z_k\log z_k.
\end{equation}
Note however that, because we do not know that $\varrho_+$ is an integrable function (only $\sqrt{\alpha}\varrho_+$ is due to \eqref{eq:phaseWeighted_1}), we only have
\[
z_k\rightharpoonup R+(1-\alpha)\varrho_+
=\alpha\varrho_++(1-\alpha)\varrho_+
=\varrho_+,
\]
on some subset $E_j=\{\varrho_+\leq j\}$.
Now, passing to the limit $k\to\infty$ in \eqref{loglog_2}, we obtain
\[
\begin{aligned}
z_k\log z_k
&\rightharpoonup
\overline{r\log r}+(1-\alpha)\varrho_+\log\varrho_+
=\varrho_+\log\varrho_+,
\end{aligned}
\]
weakly in $L^1(E_j)$. Because $z\log z$ is a convex, using \cite[Theorem~11.27]{FeireislNovotny2017}, on each
$E_j=\{\varrho_+\leq j\}$, we therefore obtain
$z_k\to \varrho_+$ a.e. on $E_j\cap\{\varrho_+>0\}$.
On $\{\varrho_+=0\}$, nonnegativity and weak $L^1$ convergence
give $\int_{\{\varrho_+=0\}}z_k\dxt\to0$, hence an almost-everywhere
convergent subsequence there as well. From here, using the diagonal argument we obtain
\begin{equation}\label{zrho}
z_k\longrightarrow\varrho_+
\quad\text{a.e. on }Q_T.
\end{equation}
The final argument requires to see that for every $b\in C_b([0,\infty))$, we have
\begin{equation}\label{earlier}
 s_k|b(\varrho_k)-b(\varrho_+)|
 =|b(z_k)-b(\varrho_+)|.
\end{equation}
It follows from the definition of \(z_k\): if \(s_k=1\), then \(z_k=\varrho_k\), and both sides equal \(|b(\varrho_k)-b(\varrho_+)|\); if \(s_k=0\), then \(z_k=\varrho_+\), so both sides are zero.

On the other hand, the convergence \eqref{zrho} and boundedness of $b$ imply that
\[
|b(z_k)-b(\varrho_+)|\longrightarrow0
 \qquad\text{in }L^1(Q_T).
\]
Combining this with \eqref{earlier} and $s_k\weakstar\alpha$ proves the first identity in
\eqref{eq:mainatoms},  as a weak-* limit in $L^\infty(Q_T)$. 
The other equality follows the same steps repeated  for $1-s_k,q_k,Q,\varrho_-$ instead of $s_k,r_k,R,\varrho_+$.

\end{proof}

\subsection{The limiting pressure and the equation for the fraction}\label{sec:gap}
We first identify the pressure $\overline{p}$. From \eqref{lim_pres} We already know that
\(
p_+(r_k)\rightharpoonup \overline{p_+(r)}.
\)
We want to prove that \(\overline{p_+(r)}=\alpha p_+(\varrho_+)\).
Because Proposition \ref{prop:phasemarginal} works only for the bounded continuous function $b$, we will use it for 
 $b_L(z)=\min\{p_+(z),L\}$.

Since $s_k\in\{0,1\}$,
$r_k=s_k\varrho_k$ and $p_+(0)=0$, we have
\[
 s_k b_L(\varrho_k)=b_L(r_k).
\]
Thus \eqref{eq:mainatoms} gives, for fixed $L$,
\[
 b_L(r_k)
 \weakstar\alpha b_L(\varrho_+)
 \quad\text{in }L^\infty(Q_T).
\]

We next estimate the error introduced by truncation. The difference
$p_+(r_k)-b_L(r_k)$ vanishes where $p_+(r_k)\leq L$.
On the remaining set, $p_+(r_k)\leq L^{-\delta}p_+(r_k)^{1+\delta}$,
so
\[
 0\leq p_+(r_k)-b_L(r_k)\leq L^{-\delta}p_+(r_k)^{1+\delta}
 \leq L^{-\delta}p_k^{1+\delta}.
\]
Using \eqref{eq:extra}, we therefore obtain
\[
 0\leq\int_{Q_T}
 \bigl[p_+(r_k)-b_L(r_k)\bigr]\dxt
 \leq C_T L^{-\delta},
\]
uniformly in $k$. In particular, for every bounded nonnegative test function $\phi$,
\[
 \int_{Q_T}\phi
   [\overline{p_+(r)}-\alpha b_L(\varrho_+)]\dxt
 =\lim_k\int_{Q_T}\phi[p_+(r_k)-b_L(r_k)]\dxt\geq0.
\]
Hence $0\leq\alpha b_L(\varrho_+)\leq\overline{p_+(r)}$
almost everywhere. 
Now, for  $L\to\infty$ functions
$\alpha\min\{p_+(\varrho_+),L\}$ increase to
$\alpha p_+(\varrho_+)$ and are uniformly bounded from above by $\overline{p_+(r)}\in L^1(Q_T)$.
Thus, by the monotone convergence theorem one gets
\[
 \overline{p_+(r)}-\alpha p_+(\varrho_+)\geq0,\qquad
 \int_{Q_T}\bigl[\overline{p_+(r)}-\alpha p_+(\varrho_+)\bigr]\dxt=0,
\]
and so $\overline{p_+(r)}=\alpha p_+(\varrho_+)$ a.e. The same argument with
$b_L(z)=\min\{p_-(z),L\}$ gives the second identity below:
\begin{equation}\label{P_ident}
 \overline{p_+(r)}=\alpha p_+(\varrho_+),\qquad
 \overline{p_-(q)}=(1-\alpha)p_-(\varrho_-).
\end{equation}
Substituting this to \eqref{eq:limmomentum} concludes the limit passage in the momentum equation.

Next, we identify the right hand side in the limiting equation for $\alpha$. Taking $b=1$ in \eqref{eq:mainflux} and using $s_kp_k=p_+(r_k)$ gives
\[\overline{s F}=
 (2\mu+\lambda)\overline{s\diver u}-\overline{p_+(r)}
 =\alpha F
 =\alpha\bigl[(2\mu+\lambda)\diver U-\overline p\bigr].
\]
Substituting $\overline p=\overline{p_+(r)}+\overline{p_-(q)}$ and the identified pressures \eqref{P_ident},
we obtain
\begin{equation}\label{wBN}
 \begin{aligned}
 (2\mu+\lambda)(\overline{s\diver u}-\alpha\diver U)
 &=\overline{p_+(r)}-\alpha \overline{p}\\
 &=(1-\alpha)\overline{p_+(r)}-\alpha \overline{p_-(q)}\\
 &=\alpha(1-\alpha)
       \bigl[p_+(\varrho_+)-p_-(\varrho_-)\bigr].
 \end{aligned}
\end{equation}

Thus $\overline{s\diver u}-\alpha\diver U=\omega_{BN}(\alpha,R,Q)\in L^2(Q_T)$, defined in \eqref{eq:BNsource}. Substituting it into  in and
\eqref{eq:limRQa} gives the balances in \eqref{eq:macro}. Substituting
$G=\alpha\diver U+\omega_{BN}(\alpha,R,Q)$ into the renormalized fraction equation  then gives \eqref{eq:BNfractionweak}.

\section{Passage to the limit in the energy inequality}\label{sec:macroenergy}
It remains to pass to the limit in \eqref{eq:wenenergy}, i.e
\begin{equation}\label{eq:wenenergy_k}
  E(r_k,q_k,u_k)(\tau)
  +\int_0^\tau\!\int_\Omega\Sstress(Du_k):Du_k\dxt\leq E^0_k,
\end{equation}
where $E^0_k$ is as given in \eqref{eq:genericinitialenergy}:
\begin{equation}
 E^0_k:=\int_\Omega
 \left[\frac{|m^0_k|^2}{2(r^0_k+q^0_k)}
                  +H_+(r^0_k)+H_-(q^0_k)\right]\dx<\infty,
\end{equation}
and to justify \eqref{eq:macroenergy}. In particular, we need to justify the appearance of a new pressure-gap dissipation term $\Rbn(\alpha,R,Q)$ on the l.h.s. of the energy inequality.

Let $0\leq\eta\in C_c^\infty(0,T)$. Multiply
\eqref{eq:wenenergy_k} by $\eta(\tau)$ and integrate in $\tau$.
Fubini's theorem gives
\begin{equation}\label{weak_ener_pass}
 \begin{aligned}
 &\int_0^T\eta(\tau)E(r_k,q_k,u_k)(\tau)\,\mathrm d\tau\\
 &\quad+\int_{Q_T}\left(\int_t^T\eta(\tau)\,\mathrm d\tau\right)
       \Sstress(Du_k):Du_k\dxt
 \leq E_k^0\int_0^T\eta(\tau)\,\mathrm d\tau.
 \end{aligned}
\end{equation}
From the discussion in Section \ref{sec:data}, we  deduced that $E^0_k\to E^0_{BN}$.
Next, using the last convergence of \eqref{eq:convectiveclosure} and the fact that
$H_\pm=p_\pm/(\gamma^\pm-1)$, we obtain
\begin{equation}\label{conv_E_hat}
 \widehat E(r_k,q_k,u_k)\rightharpoonup
 \tfrac12(R+Q)|U|^2+H_{BN}(\alpha,R,Q)\quad\text{in }L^1(Q_T).
\end{equation}

To recover the pressure-gap dissipation, note that integration by parts gives
\begin{equation}\label{div_curl}
 \int_\Omega\Sstress(Du_k):Du_k\dx
 =\int_\Omega\bigl[\mu|\curl u_k|^2+(2\mu+\lambda)|\diver u_k|^2\bigr]\dx.
\end{equation}
We again use  \cite[Theorem~11.27]{FeireislNovotny2017},
as in \eqref{eq:phaseWeighted_1} with the convex function
\[
 (a,x,y)\longmapsto\frac{x^2}{a}+\frac{y^2}{1-a},
\]
with similar lower-semicontinuous extension as used in
\eqref{eq:phaseWeighted_1}: with value $0$ when both numerator and
denominator are zero, and $+\infty$ when only the denominator is zero.
The weak convergences give
\[
 \begin{aligned}
 &(s_k,s_k\diver u_k,
       (1-s_k)\diver u_k)\rightharpoonup
 (\alpha,\overline{s\diver u},(\diver U-\overline{s\diver u}))
 \quad\text{in }L^2(Q_T;\R^3).
 \end{aligned}
\]
Since $s_k\in\{0,1\}$,
\[
 \frac{|s_k\diver u_k|^2}{s_k}
 +\frac{|(1-s_k)\diver u_k|^2}{1-s_k}
 =|\diver u_k|^2.
\]
The lower-semicontinuity therefore yields
\[
 \liminf_k\int_{Q_T}|\diver u_k|^2\dxt
 \geq\int_{Q_T}\left[
 \frac{\overline{s\diver u}^2}{\alpha}
 +\frac{(\diver U-\overline{s\diver u})^2}{1-\alpha}\right]\dxt.
\]
In particular, the limiting integral is finite, thus $\overline{s\diver u}=0$ on $\{\alpha=0\}$ and $\overline{s\diver u}=\diver U$ on
$\{\alpha=1\}$. For the remaining values, $0<\alpha<1$, we first check the algebraic identity
\[ \frac{\overline{s\diver u}^2}{\alpha}+\frac{(\diver U-\overline{s\diver u})^2}{1-\alpha}
 =|\diver U|^2
   +\frac{(\overline{s\diver u}-\alpha\diver U)^2}{\alpha(1-\alpha)}.\]
Using the identity \eqref{wBN}, we further obtain
\[
 \begin{aligned}
\frac{(\overline{s\diver u}-\alpha\diver U)^2}{\alpha(1-\alpha)}&=\frac{\alpha(1-\alpha)}{(2\mu+\lambda)^2}
     [p_+(\varrho_+)-p_-(\varrho_-)]^2\\
 &=\frac{\Rbn(\alpha,R,Q)}{2\mu+\lambda},
 \end{aligned}
\]
with $\Rbn$ defined in \eqref{eq:BNrelax}. We therefore obtained that 
\[\liminf_k\int_{Q_T}(2\mu+\lambda)|\diver u_k|^2\dxt
 \geq\int_{Q_T}[(2\mu+\lambda)|\diver U|^2+\Rbn(\alpha,R,Q)]\dxt\]
Adding the weak lower-semicontinuity bound for
$\mu\int_{Q_T}|\curl u_k|^2\dxt$, and using the identity \eqref{div_curl}, we therefore obtain

\begin{equation}\label{eq:dissipationlsc}
 \liminf_k\int_{Q_T}\Sstress(Du_k):Du_k\dxt
 \geq\int_{Q_T}\bigl[\Sstress(DU):DU+\Rbn(\alpha,R,Q)\bigr]\dxt.
\end{equation}
In particular, $\Rbn(\alpha,R,Q)\in L^1(Q_T)$.

We conclude by using the weak space--time convergence of the energy density \eqref{conv_E_hat}, the
weighted version of \eqref{eq:dissipationlsc} (with the weight $\int_t^T\eta(\tau)\,\mathrm d\tau$ as in \eqref{weak_ener_pass}), and
$E_k^0\to E_{BN}^0$. Applying Fubini's theorem once more yields
\[
 \int_0^T\eta(\tau)\left[
 E_{BN}(\alpha,R,Q,U)(\tau)
 +\int_0^\tau\!\int_\Omega
    [\Sstress(DU):DU+\Rbn(\alpha,R,Q)]\dxt
 -E_{BN}^0\right]\,\mathrm d\tau\leq0.
\]
Since $\eta\geq0$ is arbitrary, \eqref{eq:macroenergy} holds for
almost every $\tau\in(0,T)$.

This completes the proof of
Theorem~\ref{thm:main}. \hfill$\square$

\section*{Acknowledgments}
E.Z. acknowledges support
from the EPSRC Early Career Fellowship EP/V000586/1.

\section*{Declaration of AI use}

We used OpenAI’s ChatGPT in developing and writing parts of the manuscript. Starting from our drafts and in response to our prompts, ChatGPT proposed a proof of key Proposition 4.1 and worked out intermediate steps. We developed the final proof through a series of exchanges with the model. The authors take responsibility for the mathematical arguments and the final text.

\bibliographystyle{apalike}
\bibliography{Biblio1_cleaned2}

@InCollection{BurteaGavrilyukPerrin2024,
  author        = {Burtea, Cosmin and Gavrilyuk, Sergey L. and Perrin, Charlotte},
  title         = {{H}amilton's principle of stationary action in multiphase flow modeling},
  booktitle     = {Mathematical methods and modeling for mixtures of fluids and interface evolution},
  publisher     = {Paris: Soci{\'e}t{\'e} Math{\'e}matique de France (SMF)},
  year          = {2024},
  pages         = {77--106},
  isbn          = {978-2-85629-987-6},
  language      = {English},
  zbl           = {1557.76095},
  zbmath        = {7958501},
}

@Article{VasseurWenYu2019,
  author        = {Vasseur, Alexis and Wen, Huanyao and Yu, Cheng},
  title         = {{G}lobal weak solution to the viscous two-fluid model with finite energy},
  journal       = {Journal de Mathématiques Pures et Appliquées},
  publisher     = {Elsevier BV},
  year          = {2019},
  month         = May,
  volume        = {125},
  pages         = {247--282},
  doi           = {10.1016/j.matpur.2018.06.019},
  issn          = {0021-7824},
}

@Article{NovotnyPokorny2020,
  author        = {Novotn{\`y}, Anton{\'i}n and Pokorn{\`y}, Milan},
  title         = {{W}eak solutions for some compressible multicomponent fluid models},
  journal       = {Archive for Rational Mechanics and Analysis},
  publisher     = {Springer},
  year          = {2020},
  volume        = {235},
  number        = {1},
  pages         = {355--403},
}

@Article{Serre1991,
  author        = {Serre, Denis},
  title         = {{V}ariations de grande amplitude pour la densite d'un fluide visqueux compressible},
  journal       = {Physica D: Nonlinear Phenomena},
  publisher     = {Elsevier},
  year          = {1991},
  volume        = {48},
  number        = {1},
  pages         = {113--128},
  doi           = {10.1016/0167-2789(91)90055-e},
}

@Article{BaerNunziato1986,
  author        = {Baer, MR and Nunziato, JW},
  title         = {{A} two-phase mixture theory for the deflagration-to-detonation transition {D}{D}{T} in reactive granular materials},
  journal       = {International Journal of Multiphase Flow},
  publisher     = {Elsevier},
  year          = {1986},
  month         = Nov,
  volume        = {12},
  number        = {6},
  pages         = {861--889},
  doi           = {10.1016/0301-9322(86)90033-9},
  issn          = {0301-9322},
}

@Article{BreschHuang2011,
  author        = {Bresch, Didier and Huang, X.},
  title         = {{A} {M}ulti-{F}luid {C}ompressible {S}ystem as the {L}imit of {W}eak {S}olutions of the {I}sentropic {C}ompressible {N}avier–{S}tokes {E}quations},
  journal       = {Archive for Rational Mechanics and Analysis},
  publisher     = {Springer Science and Business Media LLC},
  year          = {2011},
  month         = Feb,
  volume        = {201},
  number        = {2},
  pages         = {647--680},
  doi           = {10.1007/s00205-011-0400-8},
  issn          = {1432-0673},
}

@Article{BreschHillairet2015,
  author        = {Bresch, D. and Hillairet, M.},
  title         = {{N}ote on the derivation of multi-component flow systems},
  journal       = {Proceedings of the American Mathematical Society},
  publisher     = {American Mathematical Society (AMS)},
  year          = {2015},
  month         = Apr,
  volume        = {143},
  number        = {8},
  pages         = {3429--3443},
  doi           = {10.1090/proc/12614},
  issn          = {1088-6826},
}

@Article{BreschHillairet2019,
  author        = {Bresch, Didier and Hillairet, Matthieu},
  title         = {{A} compressible multifluid system with new physical relaxation terms},
  journal       = {Annales scientifiques de l’École normale supérieure},
  publisher     = {Societe Mathematique de France},
  year          = {2019},
  volume        = {52},
  number        = {2},
  pages         = {255--295},
  doi           = {10.24033/asens.2387},
  issn          = {1873-2151},
}

@Book{DrewPassman1999,
  author        = {Drew, Donald A and Passman, Stephen L},
  title         = {{T}heory of {M}ulticomponent {F}luids},
  journal       = {Applied Mathematical Sciences},
  publisher     = {Springer New York},
  year          = {1999},
  doi           = {10.1007/b97678},
  isbn          = {9780387226378},
  issn          = {2196-968X},
}

@InBook{Gavrilyuk2011,
  author        = {Gavrilyuk, Sergey},
  title         = {{M}ultiphase flow modeling via {H}amilton's principle},
  booktitle     = {Variational models and methods in solid and fluid mechanics},
  publisher     = {Springer},
  year          = {2011},
  pages         = {163--210},
  doi           = {10.1007/978-3-7091-0983-0_4},
  isbn          = {9783709109830},
  issn          = {0254-1971},
}

@Article{GavrilyukSaurel2002,
  author        = {Gavrilyuk, Sergey and Saurel, Richard},
  title         = {{M}athematical and {N}umerical {M}odeling of {T}wo-{P}hase {C}ompressible {F}lows with {M}icro-{I}nertia},
  journal       = {Journal of Computational Physics},
  publisher     = {Elsevier},
  year          = {2002},
  month         = Jan,
  volume        = {175},
  number        = {1},
  pages         = {326--360},
  doi           = {10.1006/jcph.2001.6951},
  issn          = {0021-9991},
}

@Book{IshiiHibiki2011,
  author        = {Ishii, Mamoru and Hibiki, Takashi},
  title         = {{T}hermo-{F}luid {D}ynamics of {T}wo-{P}hase {F}low},
  publisher     = {Springer New York},
  year          = {2011},
  doi           = {10.1007/978-1-4419-7985-8},
  isbn          = {9781441979858},
}

@Article{SaurelGavrilyukRenaud2003,
  author        = {Saurel, Richard and Gavrilyuk, Sergey and Renaud, Fran{\c{c}}ois},
  title         = {{A} multiphase model with internal degrees of freedom: application to shock–bubble interaction},
  journal       = {Journal of Fluid Mechanics},
  publisher     = {Cambridge University Press},
  year          = {2003},
  month         = Nov,
  volume        = {495},
  pages         = {283--321},
  doi           = {10.1017/s002211200300630x},
  issn          = {1469-7645},
}

@Article{Novotny2020,
  author        = {Novotn{\`y}, Antonin},
  title         = {{W}eak solutions for a bi-fluid model for a mixture of two compressible non interacting fluids},
  journal       = {Science China Mathematics},
  publisher     = {Springer},
  year          = {2020},
  volume        = {63},
  number        = {12},
  pages         = {2399--2414},
}

@Article{BreschBurteaLagoutiere2022,
  author        = {Bresch, Didier and Burtea, Cosmin and Lagoutière, Frédéric},
  title         = {{M}athematical justification of a compressible bifluid system with different pressure laws: a continuous approach},
  journal       = {Applicable Analysis},
  publisher     = {Taylor \& Francis},
  year          = {2022},
  month         = July,
  volume        = {101},
  number        = {12},
  pages         = {4235--4266},
  doi           = {10.1080/00036811.2022.2103679},
  issn          = {1563-504X},
}

@Article{BreschBurteaLagoutiere2023,
  author        = {Didier Bresch and Cosmin Burtea and Frédéric Lagoutière},
  title         = {{M}athematical justification of a compressible bi-fluid system with different pressure laws: {A} semi-discrete approach and numerical illustrations},
  journal       = {Journal of Computational Physics},
  publisher     = {Elsevier BV},
  year          = {2023},
  month         = Oct,
  volume        = {490},
  pages         = {112259},
  doi           = {10.1016/j.jcp.2023.112259},
  url           = {https://www.sciencedirect.com/science/article/pii/S0021999123003546},
  issn          = {0021-9991},
}

@Article{Hillairet2007b,
  author        = {Hillairet, M},
  title         = {{P}ropagation of density-oscillations in solutions to the barotropic compressible {N}avier-{S}tokes system},
  journal       = {Journal of Mathematical Fluid Mechanics},
  publisher     = {Springer},
  year          = {2007},
  volume        = {9},
  number        = {3},
  pages         = {343--376},
}

@Book{Lions1998,
  author        = {Lions, Pierre-Louis},
  title         = {{M}athematical topics in fluid mechanics. {V}ol. 2},
  series        = {Oxford Lecture Series in Mathematics and its Applications},
  publisher     = {The Clarendon Press, Oxford University Press, New York},
  year          = {1998},
  volume        = {10},
  pages         = {xiv+348},
  isbn          = {0-19-851488-3},
  note          = {Compressible models, Oxford Science Publications},
  mrclass       = {76-02 (35-02 35Q30 76N10)},
}

@Article{Hoff1986,
  author        = {Hoff, David},
  title         = {{C}onstruction of solutions for compressible, isentropic {N}avier-{S}tokes equations in one space dimension with nonsmooth initial data},
  journal       = {Proceedings of the Royal Society of Edinburgh Section A: Mathematics},
  publisher     = {Royal Society of Edinburgh Scotland Foundation},
  year          = {1986},
  volume        = {103},
  number        = {3-4},
  pages         = {301--315},
}

@Article{Hoff1995b,
  author        = {Hoff, David},
  title         = {{S}trong convergence to global solutions for multidimensional flows of compressible, viscous fluids with polytropic equations of state and discontinuous initial data},
  journal       = {Archive for rational mechanics and analysis},
  publisher     = {Springer},
  year          = {1995},
  volume        = {132},
  number        = {1},
  pages         = {1--14},
  doi           = {10.1007/bf00390346},
}

@Article{Hoff1995a,
  author        = {Hoff, David},
  title         = {{G}lobal solutions of the {N}avier-{S}tokes equations for multidimensional compressible flow with discontinuous initial data},
  journal       = {Journal of Differential Equations},
  publisher     = {Elsevier},
  year          = {1995},
  volume        = {120},
  number        = {1},
  pages         = {215--254},
  doi           = {10.1006/jdeq.1995.1111},
}

@Article{Hoff2002,
  author        = {Hoff, David},
  title         = {{D}ynamics of singularity surfaces for compressible, viscous flows in two space dimensions},
  journal       = {Communications on Pure and Applied Mathematics: A Journal Issued by the Courant Institute of Mathematical Sciences},
  publisher     = {Wiley Online Library},
  year          = {2002},
  volume        = {55},
  number        = {11},
  pages         = {1365--1407},
  doi           = {10.1002/cpa.10046},
}

@Book{PlotnikovSokoowski2012,
  author        = {Plotnikov, Pavel and Soko{\l}owski, Jan},
  title         = {{C}ompressible {N}avier-{S}tokes equations: theory and shape optimization},
  publisher     = {Springer Science \& Business Media},
  year          = {2012},
  volume        = {73},
}

@Article{AmosovZlotnik1996a,
  author        = {Amosov, A. A. and Zlotnik, A. A.},
  title         = {{A}n estimate of the error of quasi-averaging of the equations of motion of a viscous barotropic medium with rapidly oscillating data},
  journal       = {Comput. Math. Math. Phys.},
  year          = {1996},
  volume        = {36},
  number        = {10},
  pages         = {1415--1428},
  issn          = {0965-5425},
  fjournal      = {Computational Mathematics and Mathematical Physics},
  language      = {English},
  zbl           = {0923.76263},
  zbmath        = {1205670},
}

@Article{AmosovZlotnik1996b,
  author        = {Amosov, A. A. and Zlotnik, A. A.},
  title         = {{O}n quasi-averaged equations of the one-dimensional motion of a viscous barotropic medium with rapidly oscillating data},
  journal       = {Comput. Math. Math. Phys.},
  year          = {1996},
  volume        = {36},
  number        = {2},
  pages         = {203--220},
  issn          = {0965-5425},
  fjournal      = {Computational Mathematics and Mathematical Physics},
  language      = {English},
  zbl           = {1027.35504},
  zbmath        = {2025718},
}

@Article{PerrierGutierrez2021,
  author        = {Perrier, Vincent and Gutiérrez, Enrique},
  title         = {{D}erivation and {C}losure of {B}aer and {N}unziato {T}ype {M}ultiphase {M}odels by {A}veraging a {S}imple {S}tochastic {M}odel},
  journal       = {Multiscale Modeling and Simulation},
  publisher     = {Society for Industrial & Applied Mathematics (SIAM)},
  year          = {2021},
  month         = jan,
  volume        = {19},
  number        = {1},
  pages         = {401--439},
  doi           = {10.1137/19m1306609},
  issn          = {1540-3467},
}

@article{Zlotnik2027,
title = {On weak solutions to the 1d compressible Navier-Stokes equations: A Lipschitz continuous dependence on data in weaker norms and an error of their homogenization},
journal = {Discrete and Continuous Dynamical Systems - S},
volume = {33},
number = {0},
pages = {289-321},
year = {2027},
issn = {1937-1632},
doi = {10.3934/dcdss.2026166},
url = {https://www.aimsciences.org/article/id/6a2fc5afd97cf26b6ea6c7a2},
author = {Alexander Zlotnik}
}

@Article{BakhvalovEglit1983,
  author        = {Bakhvalov, N. S. and \`Eglit, M. \`E.},
  title         = {{P}rocesses in periodic media not describable in terms of mean characteristics},
  journal       = {Dokl. Akad. Nauk SSSR},
  year          = {1983},
  volume        = {268},
  number        = {4},
  pages         = {836--840},
  issn          = {0002-3264},
  fjournal      = {Doklady Akademii Nauk SSSR},
  mrclass       = {82A55 (73E99)},
  mrnumber      = {693217},
}

@Article{DiPernaLions1989,
  author        = {DiPerna, R. J. and Lions, P. L.},
  title         = {{O}rdinary differential equations, transport theory and {S}obolev spaces},
  journal       = {Inventiones Mathematicae},
  publisher     = {Springer Science and Business Media LLC},
  year          = {1989},
  month         = Oct,
  volume        = {98},
  number        = {3},
  pages         = {511--547},
  doi           = {10.1007/bf01393835},
  issn          = {1432-1297},
}

@Article{Bresch_et_al2026,
  author        = {Bresch, Didier and Burtea, Cosmin and Gonin--Joubert, Pierre and Lagoutière, Frédéric},
  title         = {{J}ustification of a two-phase averaged system for ideal gases without heat conduction},
  journal       = {Annales de l’Institut Henri Poincaré C, Analyse non linéaire},
  publisher     = {European Mathematical Society - EMS - Publishing House GmbH},
  year          = {2026},
  month         = May,
  doi           = {10.4171/aihpc/179},
  url           = {https://arxiv.org/abs/2407.16720},
  issn          = {1873-1430},
  primaryclass  = {math.AP},
}

@Book{NovotnyStraskraba2004,
  author        = {Novotn\'y, A. and Straskraba, I.},
  title         = {{I}ntroduction to the mathematical theory of compressible flow},
  publisher     = {Oxford University Press on Demand},
  year          = {2004},
}

@PhdThesis{Lagoutiere2000,
  author = {Lagouti{\`e}re, Fr{\'e}d{\'e}ric},
  school = {Universit{\'e} Pierre et Marie Curie},
  title  = {{M}od{\'e}lisation math{\'e}matique et r{\'e}solution num{\'e}rique de probl{\`e}mes de fluides compressibles {\`a} plusieurs constituants},
  year   = {2000},
  type   = {Th{\`e}se de doctorat},
}

@Article{ChaudhuriPokornyZatorska2026,
  author  = {Chaudhuri, Nilasis and Pokorn{\`y}, Milan and Zatorska, Ewelina},
  journal = {arXiv preprint arXiv:2608.13348},
  title   = {Dissipative measure-valued solutions and weak--strong uniqueness for a viscous Baer--Nunziato system with pressure relaxation},
  year    = {2026},
}

@Article{Wen2021,
  author    = {Wen, Huanyao},
  journal   = {Calculus of Variations and Partial Differential Equations},
  title     = {On global solutions to a viscous compressible two-fluid model with unconstrained transition to single-phase flow in three dimensions},
  year      = {2021},
  number    = {4},
  pages     = {158},
  volume    = {60},
  publisher = {Springer},
}

@Book{FeireislNovotny2017,
  author    = {Feireisl, Eduard and Novotn{\'y}, Anton{\'i}n},
  publisher = {Birkh{\"a}user},
  title     = {Singular Limits in Thermodynamics of Viscous Fluids},
  year      = {2017},
  address   = {Cham},
  edition   = {2},
  isbn      = {978-3-319-63781-5},
  series    = {Advances in Mathematical Fluid Mechanics},
  doi       = {10.1007/978-3-319-63781-5},
}

@Misc{FeireislPokorny,
  author       = {Feireisl, Eduard and Pokorn{\'y}, Milan},
  howpublished = {\url{https://www2.karlin.mff.cuni.cz/~pokorny/LectureNotes/Feireisl_Pokorny_Compressible.pdf}},
  note         = {Lecture notes},
  title        = {Mathematical Theory of Compressible Viscous Fluids},
}

@Article{Maltese_et_al2016,
  author  = {Maltese, David and Mich{\'a}lek, Martin and Mucha, Piotr B. and Novotn{\'y}, Anton{\'i}n and Pokorn{\'y}, Milan and Zatorska, Ewelina},
  journal = {Journal of Differential Equations},
  title   = {Existence of Weak Solutions for Compressible Navier--Stokes Equations with Entropy Transport},
  year    = {2016},
  number  = {8},
  pages   = {4448--4485},
  volume  = {261},
  doi     = {10.1016/j.jde.2016.06.029},
}

@Article{NovotnyPokorny2021,
  author  = {Novotn{\'y}, Anton{\'i}n and Pokorn{\'y}, Milan},
  journal = {Journal of Evolution Equations},
  title   = {Continuity Equation and Vacuum Regions in Compressible Flows},
  year    = {2021},
  number  = {3},
  pages   = {2891--2922},
  volume  = {21},
  doi     = {10.1007/s00028-021-00704-3},
}
\end{document}